\documentstyle{article}
\title{Noncommutative spectral geometry of
Riemannian foliations}
\author{Yuri A. Kordyukov\\Department of Mathematics,\\
Ufa State Aviation Technical University,\\K.Marx str. 12,
Ufa 450000, Russia\\ e-mail: yurikor@math.ugatu.ac.ru}
\date{}
\def\Bbb#1{{\bf #1}}
\newtheorem{theorem}{Theorem}

\newtheorem{definition}{Definition}
\newtheorem{example}{Example}
\newtheorem{lemma}{Lemma}
\newtheorem{proposition}{Proposition}
\newtheorem{remark}{Remark}
\newenvironment{proof}{\begin{trivlist}{\item[] \bf Proof.\/}}
{\end{trivlist}}
\newtheorem{acknowledgement}{Acknowledgement}

\begin{document}
\maketitle
\newcommand{\tr}{{\rm tr}\;}
\newcommand{\trG}{{\rm tr}_G\;}
\newcommand{\trF}{{\rm tr}_{\cal F}\;}
\newcommand{\Tr}{{\rm Tr}\;}
\newcommand{\TR}{{\rm TR}\;}
\newcommand{\res}{{\rm res}\;}

\setcounter{section}{-1}
\section{Introduction}
\label{noncom}
According to \cite{Co-M,spview}, the initial datum of noncommutative
differential geometry is a spectral
triple $({\cal A},{\cal H},D)$ (see Section~\ref{proof} for
the definition), which provides a description of the corresponding
geometrical space in terms of spectral data of geometrical operators
on this space.

The purpose of this paper is to construct spectral triples given
by transversally elliptic operators with respect to a foliation on
a compact manifold and describe its dimension. The first result of the
paper is the following theorem:
\begin{theorem}
\label{triple}
Given a closed foliated manifold $(M,{\cal F})$, let a triple 
$(A,{\cal H},D)$ be defined as follows:
\begin{enumerate}
\item ${\cal A}$ is the involutive algebra $C^{\infty}_c(G_{\cal F})$ 
of smooth, compactly supported functions on the holonomy groupoid
$G_{\cal F}$ of the foliation ${\cal F}$;
\item ${\cal H}$
is the Hilbert space $L^2(M,E)$ of $L^2$-sections of a holonomy
equivariant Hermitian vector bundle $E$ equipped with the
$\ast$-representation $R_E$ of the algebra ${\cal A}$ (\ref{K});
\item $D$ is a first order self-adjoint transversally
elliptic operator in $L^2(M,E)$ with the holonomy invariant
transversal principal symbol such that the operator $D^2$ is
self-adjoint and has the scalar principal symbol.
\end{enumerate}
Then $({\cal A},{\cal H},D)$ is a finite-dimensional spectral triple.
\end{theorem}

A geometrical example of spectral triples considered in Theorem~\ref{triple}
is given by the transverse signature operator on a Riemannian foliation.

\begin{example} 
Let $(M,{\cal F})$ be a Riemannian foliation,
equipped with a bundle-like metric $g_M$. Let $F=T{\cal F}$ be the
integrable distribution in $TM$ of tangent $p$-planes to the foliation,
and $H=F^{\bot}$ be the orthogonal complement to $F$. So we have
a decomposition of $TM$ into a direct sum $TM=F\oplus H$ and
the corresponding decomposition
of the de Rham differential $d$ in the form
$d=d_F+d_H+\theta$, 
where the tangential de Rham differential $d_F$
and the transversal de Rham differential $d_H$ are first order differential
operators, and $\theta$ is zeroth order.

To define a spectral triple $({\cal A},{\cal H},D)$,
we take the Hilbert space ${\cal H}$
to be the space $L^2(M,\Lambda^{*}H^{*})$ of transversal
differential forms, equipped with the natural action of the
algebra ${\cal A}=C^{\infty}_c(G)$, and the operator $D$ to be
the transverse signature operator $d_H+\delta_H$ (see Section~\ref{self}
for more details).
\end{example}

Transversally elliptic operators on manifolds, equipped with an action
of a compact Lie group, were introduced by Atiyah and Singer in \cite{Atiyah}.
In the context of noncommutative differential geometry, these operators
appeared in \cite{CoNG} to provide examples of Fredholm modules, associated
with foliated manifolds. Namely, it was proved there that any zeroth order
transversally elliptic operator with the holonomy invariant transversal
principal symbol gives rise to a finite-dimensional Fredholm module over
the foliation algebra $C^{\infty}_c(G_{\cal F})$ (see also \cite{Co,H-Sk}).
Theorem~\ref{triple} provides an extension of the above mentioned
result to the case of transversally elliptic operators of positive order.

The next problem is to describe dimension of the spectral triples in question.
The usual notion of dimension for a general spectral triple
$({\cal A},{\cal H},D)$ (\cite{Co}) is given by the degree of summability $d$
of the operator $(D-i)^{-1}$, that is, by
the least $p$ such that the operator $a(D-i)^{-1}, a\in {\cal A}$
is an operator of the Schatten ideal ${\cal L}^p({\cal H})$.
In the case under consideration, $d$ is equal to the codimension
$q$ of the foliation ${\cal F}$ (see Proposition~\ref{invers} for a
proof). If we are looking at a geometrical space as a union of pieces
of different dimensions, this notion of dimension of the
corresponding spectral triple gives only an upper bound on 
dimensions of various pieces. To take into account
lower dimensional pieces of the space under consideration,
Connes and Moscovici \cite{Co-M} suggested
that the correct notion of dimension is given not by
a single real number $d$ but by a subset ${\rm Sd}\subset
{\Bbb C}$, which is called the dimension spectrum of
the given triple (see Section~\ref{proof} for the definition).
\medskip
\par
The second result of the paper is a description of the spectrum dimension
of the spectral triples defined in Theorem~\ref{triple}.

\begin{theorem}
\label{dim}
A spectral triple $({\cal A},{\cal H},D)$ as in Theorem~\ref{triple} 
has discrete dimension spectrum ${\rm Sd}$, which is
contained in the set $\{v\in {\Bbb N}:v\leq q\}$ and is simple.
\end{theorem}

In \cite{trans}, the author studied analytic properties of transversally
elliptic operators with respect to noncompact Lie group actions.
In particular, results of \cite{trans} allows us to define
(finite-dimensional) Fredholm modules given by transversally elliptic
operators with the invariant transversal principal symbol on a smooth
manifold equipped with a Lie group action and claim that the
corresponding spectral triples have discrete dimension spectrum
(but the dimension spectrum might be not simple, if there are singular
orbits).
In this paper, we combine general methods
of \cite{trans} with a further elaboration of pseudodifferential calculus
on foliated manifolds \cite{asymp,tang} to give a more precise description of
the dimension spectrum for spectral triples associated with foliated
manifolds.

This work concerns to the simplest examples of spectral triples associated
with foliated manifolds.
In the forthcoming paper \cite{noncom2}, we will extend our considerations
to foliated manifolds, equipped with a triangular transversal structure
as in \cite{Co-M}, using a transversal
pseudodifferential calculus, modelled on the Beals-Greiner pseudodifferential
calculus on Heisenberg manifolds (see \cite{B-G,Co-M}).
\medskip\par
Contents:
\medskip\par
1. Transversal pseudodifferential calculus

\hskip 12pt 1.1. Preliminaries

\hskip 12pt 1.2. Classes $\Psi^{m,-\infty}(M,{\cal F},E)$

\hskip 12pt 1.3. Anisotropic Sobolev spaces and classes
$\Psi^{m,\mu}(M,{\cal F},E)$

\hskip 12pt 1.4. Symbolic properties of $\Psi^{m,-\infty}(M,{\cal F},E)$

\hskip 12pt 1.5. Residue trace

2. Transversally elliptic operators

\hskip 12pt 2.1. Definition and basic properties

\hskip 12pt 2.2. Complex powers

\hskip 12pt 2.3. $G$-trace

\hskip 12pt 2.4. Zeta-function

3. Spectral triples of Riemannian foliations

\hskip 12pt 3.1. Proof of main theorems

\hskip 12pt 3.2. Geometric example

\hskip 12pt 3.3. Concluding remarks

\section{Transversal pseudodifferential calculus}
\subsection{Preliminaries}
Throughout in the paper, we consider a closed, connected, oriented
foliated manifold
$(M,{\cal F})$, $\dim M = n$, $\dim {\cal F} = p$, $p + q = n$, and 
a complex vector bundle
$E$ on $M$ of rank $r$. We fix  a
Riemannian metric on $M$ with the corresponding distance $\rho $ 
and an Hermitian structure on $E$.

We will denote by $G=G_{\cal F}$ the holonomy groupoid of $(M,{\cal F})$.
$G$ is equipped with the source and the target
maps $s,r:G\rightarrow M$. We will make use of standard notation:
$G^{(0)}=M$ is the set of objects, $G^x=\{\gamma\in G:r(\gamma)=x\}$,
$G^x_x=\{\gamma\in G : s(\gamma)=r(\gamma)=x\}$, $x\in M$.
For any $x\in M$, $s$ defines a covering map from $G^x$ to the leaf through
the point $x$ associated with the holonomy group $G^x_x$ of the leaf.
We will identify a point $x\in  M$ with the identity element in $G^{x}_{x}$.
Let $dx$ be the Riemannian volume form on $M$,
$\lambda _{L}$ the Riemannian volume form on a leaf $L$ of ${\cal F}$
and, for any $x\in M$, $\lambda ^{x}$ its lift to a density
on the holonomy covering $G^x$.
We will make use of notation $(x,y)\in I^p\times I^q$
($I=(-1,1)$) for the local coordinates given by a foliated chart
$\kappa: I^p\times I^q\rightarrow M$
and $(\xi,\eta) \in {\Bbb R}^p\times {\Bbb R}^q$ for the dual
coordinates (in $T^*M$).

The holonomy groupoid $G$ has a structure of a smooth manifold
of dimension $2p+q$. Recall briefly the construction
of an atlas on $G$ \cite{Co79}.
Let $\kappa: I^p\times I^q\rightarrow M, \kappa': I^p\times I^q\rightarrow M$,
be two foliated charts,
$\pi: \{0\}\times I^q\rightarrow \pi(\{0\}\times I^q)
=D^q\subset M$, $\pi': \{0\}\times I^q\rightarrow \pi'(\{0\}\times I^q)=D'^q
\subset M$
be the corresponding transversals to the foliation.
The foliation charts $\kappa$, $\kappa'$ are called
{\bf compartible}, if, for any points $m\in U=\kappa
(I^p\times I^q)$ and $m'\in U'=\kappa'(I^p\times I^q)$
such that $m=\kappa(x,y)$, $m'=\kappa'(x',y)$, there is
a leafwise path $\gamma$ from $m$ to $m'$ such that the corresponding
holonomy map $h_{\gamma}$ maps the germ $\pi_m$
of the transversal $\pi$ at the point $m$ to the germ $\pi'_{m'}$
of the transversal $\pi'$ at the point $m'$:
$h_{\gamma}\pi_m=\pi'_{m'}$.

For any pair of compartible foliation charts, $\kappa$ and $\kappa'$,
let $W(\kappa,\kappa')$ be a subset in $G_{\cal F}$:
$$
W(\kappa,\kappa')=\{\gamma\in G: s(\gamma)=m\in U,
r(\gamma)=m'\in U', h_{\gamma}\pi_m=\pi'_{m'}\},
$$
equipped with a coordinate map
\begin{equation}
\label{wchart}
\Gamma:I^p\times I^p\times I^q\rightarrow W(\kappa,\kappa'),
\end{equation}
which associates to any $(x,x',y)\in I^p\times I^p\times I^q$
the element $\gamma\in W(\kappa,\kappa')$ such that $s(\gamma)=m=\kappa(x,y)$,
$r(\gamma)=m'=\kappa'(x',y)$ and $h_{\gamma}\pi_m=\pi'_{m'}$.

As shown in \cite{Co79}, the coordinate patches
$W(\kappa,\kappa')$ form an atlas of a $2p+q$-dimensional
manifold on $G$.
\medskip
\par
Denote by $C^{\infty}_c(G,{\cal L}(E))$ the space
of smooth, compactly supported sections of the vector bundle
$(s^*E)^*\otimes r^*E$ on $G$.
Otherwise speaking, the value of $k\in C^{\infty}_c(G,{\cal L}(E))$ 
at any point
$\gamma \in G$ is a linear map
$k(\gamma):E_{s(\gamma)}\rightarrow E_{r(\gamma)}$.
Any element $k\in C^{\infty}_c(G,{\cal L}(E))$ defines an operator
$R_E(k): C^{\infty}(M,E)\rightarrow C^{\infty}(M,E)$ by the formula
\begin{equation}
\label{tang}
R_E(k)u(x)=\int_{G^x}k(\gamma)u(s(\gamma))d\lambda^x(\gamma),
u\in C^{\infty}(M,E), x\in M.
\end{equation}
which is said to be a tangential operator on $(M,{\cal F})$ 
defined by the tangential kernel $k$.

\subsection{Classes $\Psi^{m,-\infty}(M,{\cal F}, E)$}
In this section, we introduce the algebra  $\Psi^{*,-\infty}(M,{\cal F}, E)$
of transversal pseudodifferential operators on the foliated manifold
$(M,{\cal F})$, which can be considered as an analogue of the algebra of
pseudodifferential operators on a closed manifold.
This algebra can be realized as a Guillemin-Sternberg
algebra ${\cal R}_{\Sigma}$ \cite{GS79}, corresponding to 
a coisotropic conic submanifold $\Sigma$ in the punctured 
cotangent bundle $\tilde{T}^*M$ (see below), therefore, many its properties
can be deduced, just referring to the corresponding results for
these general algebras.
Here we prefer to give direct proofs (when it is possible), since
this is simpler and allows us to extend these results
to more general cases \cite{noncom2}.

Recall that a function $k\in C^{\infty}(I^{p} \times I^{p} \times I^{q} 
\times {\Bbb R}^{q},{\cal L}({\Bbb C}^r))$ belongs  to  the  class
$S^{m}(I^{p}\times I^{p}\times I^{q} \times {\Bbb R}^{q},
{\cal L}({\Bbb C}^r))$, if, for any multiindices $\alpha $ and $\beta $, 
there exists a
constant $C_{\alpha, \beta} > 0$ such that
$$
\| \partial^{\alpha}_{\eta} \partial^{\beta}_{(x,x',y)}k(x,x',y,\eta )\|
\leq C_{\alpha \beta }(1 +\vert \eta \vert )^{ m-\vert \alpha \vert },
(x,x',y)\in I^{p}\times I^p\times I^q, \eta \in {\Bbb R}^{q}.
$$
In what follows, we will consider only classical symbols.
Recall that a function $k\in C^{\infty}(I^{p}\times I^{p}\times I^{q}
\times {\Bbb R}^{q}, {\cal L}({\Bbb C}^r))$
is {\bf a classical symbol} of order $z\in {\Bbb C}$
($k\in S^{z,-\infty}(I^{p}\times I^{p}\times I^{q} \times {\Bbb R}^{n},
{\Bbb R}^{p}, {\cal L}({\Bbb C}^r))$), if
$k$ is represented as an asymptotic sum
$$
k(x,x',y,\eta)\sim \sum_{j=0}^{\infty} \theta(\eta)
k_{z-j}(x,x',y,\eta),
$$
where $k_{z-j}\in C^{\infty}(I^{p}\times I^{p}\times I^{q} \times
({\Bbb R}^{q}\backslash \{0\}), {\cal L}({\Bbb C}^r))$
is homogeneous in $\eta$ of degree $z-j$, that is,
$$
k_{z-j}(x,x',y,t\eta)=t^{z-j}k_{z-j}(x,x',y,\eta), t>0,
$$
and $\theta$ is a smooth function on ${\Bbb R}^{q}$ such that
$\theta(\eta)=0$ for $|\eta|\leq 1$,  $\theta(\eta)=1$ for $|\eta|\geq 2$.

A symbol $k \in S ^{m} (I ^{p} \times I^p\times I^q\times {\Bbb R}^{q}, 
{\cal L}({\Bbb C}^r))$ defines an operator $A: C^{\infty}_c(I^n, {\Bbb C}^r)
\rightarrow C^{\infty}(I^n, {\Bbb C}^r)$ by the formula
\begin{equation}
\label{loc}
Au(x,y)=(2\pi)^{-q} \int	e^{i(y-y')\eta}k(x,x',y,\eta) u(x',y')
\,dx'\,dy'\,d\eta,
\end{equation}
\noindent where \(u \in C^{\infty}_{c}(I^{n}, {\Bbb C}^r), x \in I^{p},
y \in I^{q}\). Denote by $\Psi^{m,-\infty}(I^n,I^p,{\Bbb C}^r)$
the class of operators of the form (\ref{loc}) with $k \in
S ^{m} (I ^{p} \times I^p\times I^q\times {\Bbb R}^{q}, 
{\cal L}({\Bbb C}^r))$ such that its Schwartz kernel is compactly supported 
in $I^n\times I^n$.

It is very useful to note that the algebra 
$\Psi^{*,-\infty}(I^n,I^p,{\Bbb C}^r)$
has structure of a
crossed product of an algebra of pseudodifferential operators on transversals
to the foliation (in $y$ variables) by the leafwise equivalence relation.
More precisely, it can be formulated as follows.
If we represent the space $L^2(I^n,{\Bbb C}^r)$ as the $L^2$ space of
$L^2(I^q, {\Bbb C}^r))$-valued functions on $I^p$,
$L^2(I^n,{\Bbb C}^r)= L^2(I^p, L^2(I^q, {\Bbb C}^r))$, then
an operator $A\in\Psi^{m,-\infty}(I^n,I^p,{\Bbb C}^r)$ can be written
in the following form:
\begin{equation}
\label{loc1}
A\bar{u}(x)=\int  A(x,x')\bar{u}(x')\,dx', x\in I^p,
\end{equation}
where $\bar{u}\in C_c^{\infty}(I^p, L^2(I^q, {\Bbb C}^r))$ such that
$\bar{u}(x)\in C_c^{\infty}(I^q, {\Bbb C}^r)$ for any $x\in I^p$, and,
for any $x\in I^p$, $x'\in I^p$, the operator $A(x,x')$ is a pseudodifferential
operator of order $m$ on $I^q$:
\begin{equation}
\label{loc2}
A(x,x')v(y)=(2\pi)^{-q} \int  e^{i(y-y')\eta}k(x,x',y,\eta) v(y')\,dy'\,d\eta,
v\in C^{\infty}_c(I^p,{\Bbb C}^r).
\end{equation}

The principal symbol of $A\in\Psi^{m,-\infty}(I^n,I^p,{\Bbb C}^r)$ is
defined to be a matrix-valued function $\sigma_A$ on 
$I^p\times I^p\times I^q\times ({\Bbb R}^q
\backslash \{0\})$ given by the formula
\begin{equation}
\label{loc:def}
\sigma_A(x,x',y,\eta)=k_m(x,x',y,\eta),
\end{equation}
where $k_m$ is the homogeneous (of degree $m$) component of the complete
symbol $k$ of the operator $A$.
We have the following properties of the principal symbols of
operators from $\Psi^{m,-\infty}(I^n,I^p,{\Bbb C}^r)$.
\begin{lemma}
\label{symbol}
(1) Given $A\in\Psi^{m_1,-\infty}(I^n,I^p,{\Bbb C}^r)$ and    
$B\in\Psi^{m_2,-\infty}(I^n,I^p,{\Bbb C}^r)$, the composition $C=AB$
belongs to the class $\Psi^{m_1+m_2,-\infty}(I^n,I^p,{\Bbb C}^r)$, and
its principal symbol $\sigma_C$ is given by the formula
$$
\sigma_C(x,x',y,\eta)=\int \sigma_A(x,x'',y,\eta)\sigma_B(x'',x',y,\eta)dx''.
$$
\medskip
\par
\noindent (2) The principal symbol $\sigma_A$ of an operator
$A\in\Psi^{m,-\infty}(I^n,I^p,{\Bbb C}^r)$ is transformed under a foliated
coordinate change $x_1=\phi(x,y), y_1=\psi(y)$, via the following
formula
$$
\sigma_{A_1}(\phi(x,y),\phi(x',y),\psi(y),(d\psi(y)^*)^{-1}(\eta))
=\sigma_A(x,x',y,\eta),
$$
where the operator $A$ is assumed to be written in the coordinates $(x,y)$
and $A_1$ denotes the operator $A$, written in the coordinates $(x_1,y_1)$.
\end{lemma}

\begin{proof}
If $A\in\Psi^{m,-\infty}(I^n,I^p,{\Bbb C}^r)$ is written in the form
(\ref{loc1}), the principal symbol $\sigma_A$ of $A$ can be expressed in
terms of the principal symbols $\sigma_{A(x,x')}$ of the operators
$A(x,x')$ as follows:
$$
\sigma_A(x,x',y,\eta)=\sigma_{A(x,x')}(y,\eta).
$$
If operators $A\in\Psi^{m_1,-\infty}(I^n,I^p,{\Bbb C}^r)$ and
$B\in\Psi^{m_2,-\infty}(I^n,I^p,{\Bbb C}^r)$ are written in the form
(\ref{loc1}) with the corresponding families
$A(x,x')\in\Psi^{m_1}(I^q,{\Bbb C}^r)$ and
$B(x,x')\in\Psi^{m_2}(I^q,{\Bbb C}^r)$ accordingly, then it is easy to
see that the composition
$C=AB$ is written in the form (\ref{loc1}) with
$$
C(x,x')=\int A(x,x'')B(x'',x')dx''.
$$
Using these facts and standard pseudodifferential calculus,
the lemma can be easily proved.
\end{proof}

If $\kappa: I^p\times I^q\rightarrow U=\kappa(I^p\times I^q)
\subset M,
\kappa': I^p\times I^q\rightarrow U'=\kappa'(I^p\times I^q)
\subset M$, are two compartible foliated charts
on $M$ equipped with trivializations of the vector bundle $E$ over them,
we can transfer an operator $A\in\Psi^{m,-\infty}(I^n,I^p,{\Bbb C}^r)$
to an operator $A': C^{\infty}_c(U,E)\rightarrow C^{\infty}_c(U',E)$, 
which extends in a trivial way to an operator in $C^{\infty}(M,E)$,
denoted also by $A'$.
The resulting operator $A'$ is said to be
{\bf an elementary operator} of class $\Psi ^{m,-\infty}(M,{\cal F},E)$.

\begin{definition} The  class $\Psi ^{m,-\infty}(M,{\cal F},E)$
consists of operators $A$, acting from $C^{\infty}(M,E)$ to $C^{\infty}(M,E)$,
such that $A$ can be represented in the form $A=\sum_{i=1}^k A_i + K$,
where $A_i$ are elementary operators
of class $\Psi ^{m,-\infty}(M,{\cal F},E)$, corresponding
to some pairs $\kappa_i,\kappa'_i$ of compartible foliated charts,
$K\in \Psi ^{-\infty}(M,E)$.
\end{definition}

To give an invariant definition of the principal symbol for operators of
$\Psi ^{m,-\infty}(M,{\cal F},E)$, let us show how these
operators can be represented as Fourier integral operators, associated with
some canonical relation on the punctured cotangent space
$\tilde{T}^*M=T^{*}M\backslash \{0\}$.

It is well-known that the foliation ${\cal F}$ can be lifted to a
foliation ${\cal F}_N$ in the punctured conormal bundle
$\tilde{N}^*{\cal F}$, which is transversally parallelizable and, therefore,
has trivial holonomy (see \cite{Molino}). In local coordinates
$(x,y,\eta)$ on $\tilde{N}^*{\cal F}$ given by a foliated chart on $M$,
plaques of the foliation ${\cal F}_N$
are defined by $y=\mbox{const},\eta=\mbox{const}$.
It is easy to see that the leaf $\tilde{L}_{\eta}$ of the
foliation ${\cal F}_N$  through a point
$\eta\in \tilde{N}^*{\cal F}$ is diffeomorphic to the holonomy
covering $G_{\cal F}^x$ of the leaf $L_x, x=\pi(\eta)$,
of the foliation ${\cal F}$ through the point $x$, therefore, we
can give the following description of the holonomy groupoid
of ${\cal F}_N$.

Recall that, for any smooth leafwise path $\gamma $ from $x\in
M$ to $y\in M$, there is defined the map 
$dh_{\gamma}^{\ast}: N^{\ast}_{y}{\cal F} 
\rightarrow N^{\ast}_{x}{\cal F}$,  being  the
codifferential of the holonomy map $h_{\gamma }$,
corresponding to $\gamma$ (cf., for instance, \cite{Co}).
The holonomy groupoid $G_{{\cal F}_N}$ of the lifted
foliation ${\cal F}_N$ consists of all
$(\gamma,\eta)\in G_{\cal F}\times \tilde{N}^*{\cal F}$ such
that $r(\gamma)=\pi(\eta)$ with the source map $s:G_{{\cal F}_N}\rightarrow
\tilde{N}^{*}{\cal F}, s(\gamma,\eta)=dh_{\gamma}^{*}(\eta)$
and the target map
$r:G_{{\cal F}_N}\rightarrow \tilde{N}^{*}{\cal F}, r(\gamma,\eta)=\eta$.
The bundle map
$\pi:\tilde{N}^{*}{\cal F}\rightarrow M$ induces a map
$\pi_G:G_{{\cal F}_N}\rightarrow G_{\cal F}$ by
$$
\pi_G(\gamma,\eta)=\gamma, (\gamma,\eta)\in
G_{{\cal F}_N}.
$$

There is also a symplectic description of the lifted foliation.
Consider $T^{*}M$ as a symplectic manifold, equipped with the canonical
symplectic structure. Then $\tilde{N}^*{\cal F}$ is a coisotropic submanifold
in $\tilde{T}^*M$, and
the foliation ${\cal F}_N$ is the corresponding null-foliation.
It is well-known that the mapping
\begin{equation}
\label{can}
(r,s):G_{{\cal F}_N}\rightarrow \tilde{T}^*M\times \tilde{T}^*M
\end{equation}
defines an immersed canonical relation in $\tilde{T}^*M$, which
is often called by the flowout of the coisotropic submanifold
$\tilde{N}^*{\cal F}$.

It can be easily checked that the algebra
of Fourier integral operators, associated with this canonical relation, is
the algebra $\Psi^{*,-\infty}(M,{\cal F},E)$ introduced above. We only need
to be more precise about the immersed canonical relation~(\ref{can}). 
Namely, let us define  the space $I^m(M\times M, G'_{{\cal F}_N})$
of compactly supported Lagrangian distributions, taking finite sums of
elementary Lagrangian distributions as we did above in the definition of
classes $\Psi^{*,-\infty}(M,{\cal F},E)$. A precise statement is that
the class $\Psi^{m,-\infty}(M,{\cal F},E)$ consists of all operators
in $C^{\infty}(M,E)$ with Schwartz kernels from the
space $I^{m-p/2}(M\times M, G'_{{\cal F}_N})$.

Now let us show how the notion of the principal symbol of operators of class
$\Psi^{m,-\infty}(M,{\cal F},E)$ as Fourier integral operators agrees
with the local definition given by (\ref{loc:def}).
According to \cite[Section 25.1]{Ho},
the principal symbol of an operator of class $\Psi^{m,-\infty}(M,{\cal F},E)$
as a Fourier integral operator is
a half-density on $G_{{\cal F}_N}$ homogeneous of degree $m+q/2$ defined
as follows.
Let $\kappa: I^p\times I^q\rightarrow U=\kappa(I^p\times I^q)
\subset M,
\kappa': I^p\times I^q\rightarrow U'=\kappa'(I^p\times I^q)
\subset M$, be two compartible foliated charts
on $M$ equipped with trivializations of the vector bundle $E$ over them.
Define a foliated coordinate map
$\Gamma_N:I^p\times I^p\times I^q\times {\Bbb R}^q\rightarrow G_{{\cal F}_N}$
by $\Gamma_N(x,x',y,\eta)=(\Gamma(x,x',y),(d\kappa'^*)^{-1}(x',y,\eta))$,
where $(x,x',y,\eta)\in I^p\times I^p\times I^q\times {\Bbb R}^q$,
$\Gamma$ is the coordinate map given by (\ref{wchart})
and $(d\kappa'^*)^{-1}:I^p\times I^q\times {\Bbb R}^q\rightarrow
\tilde{N}^*{\cal F}$ is the inverse to the codifferential of $\kappa'$.
In this coordinate chart, the half-density principal symbol
$\sigma_A$ of an elementary operator
$A\in\Psi^{m,-\infty}(M,{\cal F},E)$ (given by (\ref{loc}) in
$W(\kappa,\kappa')$) is defined to be equal to
\begin{equation}
\label{princ}
k_m(x,x',y,\eta)(dx\ dx'\ dy\ d\eta)^{1/2},
\end{equation}
where $k_m$ is the homogeneous component of the complete symbol $k$ of
degree $m$.
The half-density (\ref{princ}) can be identified with a leafwise half-density,
using the canonical transversal symplectic form $dy\ d\eta$,
and, moreover, with a smooth section from
$C^{\infty}(G_{{\cal F}_N},{\cal L}(\pi^*E))$,
using the fixed leafwise density $\lambda=\{\lambda_L : L\in M/{\cal F}\}$.

Let $S^m(G_{{\cal F}_N}, \pi^*E)$ be the space of all sections
$s\in C^{\infty}(G_{{\cal F}_N},{\cal L}(\pi^*E)), s=s(\gamma,\eta)$
homogeneous in $\eta$ of degree $m$
 such that $\pi_G({\rm supp}\,s)$ is compact in $G_{\cal F}$.
Then the principal symbol $\sigma_A$ of an operator $A\in \Psi^{m,-\infty}
(M,{\cal F},E)$ is globally defined as an element
of $S^m(G_{{\cal F}_N}, \pi^*E)$.
We can also consider the principal symbol of the operator $A$ as
the corresponding tangential operator $R_{\pi^*E}(\sigma_A)$
on $\tilde{N}^*{\cal F}$ with respect to ${\cal F}_N$.

The space
$$
S^{*}(G_{{\cal F}_N}, \pi^*E)=\bigcup_mS^m(G_{{\cal F}_N}, \pi^*E)
$$
carries a structure of an involutive algebra, defined by
its embedding into the foliation algebra
$C^{\infty}_c(G_{{\cal F}_N}, \pi^*E)$.
By Lemma~\ref{symbol}, we obtain symbolic properties of
$\Psi^{m,-\infty}(M,{\cal F},E)$ (see also \cite[Section 2]{GS79} for the
corresponding general result for Fourier integral operators).

\begin{proposition}
\label{prop}
(1) The principal symbol mapping
$$
\sigma : \Psi ^{*,-\infty}(M,{\cal F},E)\rightarrow
S^{*}(G_{{\cal F}_N}, \pi^*E)
$$
is an algebra homomorphism. Otherwise speaking,
if $A\in  \Psi ^{m_{1},-\infty}(M,{\cal F},E)$
and $B\in  \Psi ^{m_{2},-\infty}(M,{\cal F},E)$,
then \(C = AB\) belongs to
\(\Psi ^{m_{1}+m_{2},-\infty}(M,{\cal F},E)\) and
$\sigma_{AB}=\sigma_A\sigma_B$.
\medskip
\par
\noindent (2) If the density $dx$ on $M$ is holonomy invariant, then
$\sigma$ is a $*$-homomorphism of involutive algebras, i.e.,
if $A\in  \Psi ^{m,-\infty}(M,{\cal F},E)$,
then \(A^{*} \in  \Psi ^{m,-\infty}(M,{\cal F},E)\) and $\sigma_{A^*}=
(\sigma_A)^*$.
\end{proposition}

The next problem is to state $L^2$-continuity of operators
from $\Psi^{0,-\infty}(M,{\cal F},E)$. We refer the reader to
\cite[Theorem 25.3.8]{H4}
for the corresponding general result on $L^2$-continuity
of Fourier integral operators, but, indeed, our case is a model
case for this general theorem.
\begin{proposition}
Any operator $A\in \Psi^{0,-\infty}(M,{\cal F},E)$ defines a bounded operator
in the Hilbert space $L^2(M,E)$.
\end{proposition}
\begin{proof} 
The proposition follows immediately, if we make use
of a representation of the operator $A$ in the form (\ref{loc1}) and
apply the theorem on $L^2$ boundedness of zero-order
pseudodifferential operators.
\end{proof}

\subsection{Anisotropic Sobolev spaces and classes
$\Psi ^{m,\mu}(M,{\cal F},E)$}

Our norm estimates will be given in terms of the scale of
anisotropic Sobolev spaces $H^{s,k}(M,{\cal F},E)$, $s\in
{\Bbb R}, k\in {\Bbb R}$ (\cite{tang,asymp}).
Let us briefly recall its definition.

\begin{definition} The space $H^{s,k}({\Bbb R}^{n},{\Bbb R}^{p},
{\Bbb C}^r)$ consists of all $u
\in  S'({\Bbb R}^{n}, {\Bbb C}^r)$ such that its 
Fourier transform $\tilde{u} \in L^{2}_{\rm loc}({\Bbb R}^{n}, 
{\Bbb C}^r)$   and
\begin{equation}
\label{(1.9)}
\| u\|^{2}_{s,k}
= \int \int  |\tilde{u}(\xi ,\eta )|^{2}(1 +|\xi |^{2} + |\eta |^{2})^{s}
(1 +|\xi |^{2})^{k}d\xi  d\eta  < \infty.
\end{equation}
\end{definition}
The identity (\ref{(1.9)}) serves as a definition of a Hilbert
norm in $H^{s,k}({\Bbb R}^{n},{\Bbb R}^{p}, {\Bbb C}^r)$.

\begin{definition} 
The space $H^{s,k}(M,{\cal F},E)$ consists of all $u\in  {\cal D}'(M,E)$ such
that,  for any foliated coordinate chart $\kappa	: I^{p} \times I^{q}
\rightarrow  U = \kappa (I^{p} \times I^{q}) \subset M$, for any
trivialization of the bundle $E$ over it  and
for any $\phi  \in C^{\infty }_{c}(U)$,  the	function $\kappa^{\ast}(\phi u)$
belongs  to  the  space
$H^{s,k}({\Bbb R}^{n},{\Bbb R}^{p}, {\Bbb C}^r)$.
\end{definition}

Fix a finite covering $\{ U_{i} : i = 1,\ldots ,d \}$ of $M$  by
foliated coordinate patches with foliated coordinate charts $\kappa _{i} :
I^{p} \times I^{q} \rightarrow U_{i} = \kappa _{i}(I^{p} \times I^{q})$,
and a partition of unity $\{ \phi _{i} \in C^{\infty }(M):  i
= 1,\ldots,d \}$ subordinate to  this  covering.  A  scalar
product  in $H^{s,k}(M,{\cal F},E)$ can be equivalently defined by the formula
\begin{equation}
(u,v)_{s,k} =\sum ^{d}_{i=1} (\kappa^{\ast}(\phi _{i}u),
\kappa ^{\ast}(\phi _{i}v))_{s,k}, u,v\in H^{s,k}(M,{\cal F},E).
\label{(1.10)}
\end{equation}

Using the anisotropic Sobolev spaces $H^{s,k}(M,{\cal F},E),
s\in{\Bbb R}, k\in {\Bbb R}$, we can give a sufficient condition for
a bounded operator $T$ in $L^2(M,E)$ to be an operator of trace class,
adapted to the foliated structure of $M$ \cite{asymp}.
\begin{proposition}
\label{trace}
Let $T$ be a bounded operator in $L^2(M,E)$ such that
$T$ defines a bounded operator from $L^2(M,E)$
to $H^{s,k}(M,{\cal F},E)$ with some $s>q$ and
$k>p$. Then $T$ is an operator
of trace class with
the following estimate of its trace norm:
$$
\|T\|_{1}\leq C\|T:L^2(M,E)\rightarrow
H^{s,k}(M,{\cal F},E)\|.
$$
\end{proposition}

Now we introduce operator classes $\Psi ^{m,\mu}(M,{\cal F},E)$ associated
with the scale  $H^{s,k}(M,{\cal F},E)$. A local description of
these classes (denoted by $\Psi ^{m,\mu}(M,{\cal F},\delta)$ there) was 
given in \cite{asymp,tang} by means of 
H\"ormander classes of pseudodifferential operators with tempered metrics.
Here we define classes $\Psi ^{m,\mu}(M,{\cal F},E)$ globally, using
the definition of Fourier integral operators, corresponding to a pair
$(\Lambda_0,\Lambda_1)$ of
intersected Lagrangian submanifolds \cite{Gu-U,Gr-U}:
$\Lambda_0$ is the diagonal in $\tilde{T}^*M\times \tilde{T}^*M$,
$\Lambda_1$ is the holonomy groupoid $G_{{\cal F}_N}$. One of essential
advantages of our approach is possibility to make use of
the symbolic calculus for Fourier integral operators.

\begin{definition}
We say that a function $a\in C^{\infty}(I^{p}\times I^p\times I^q
 \times {\Bbb R}^{n}\times {\Bbb R}^p,{\cal L}({\Bbb C}^r)),
a=a(s,x,y,\xi,\eta,\sigma)$ belongs  to  the  class
$S^{m,\mu}(I^{n} \times {\Bbb R}^{n}, {\Bbb R}^{p},
{\cal L}({\Bbb C}^r))$, if, for any multiindices $\alpha, \beta$ and
$\gamma$, there exists a
constant $C_{\alpha, \beta, \gamma} > 0$ such that
\begin{eqnarray*}
\| \partial^{\alpha}_{(\xi,\eta)} \partial^{\beta}_{\sigma}
\partial^{\gamma}_{(s,x,y)} a(s,x,y,\xi,\eta,\sigma )\|
\leq C_{\alpha \beta \gamma}(1 +\vert \xi \vert +
\vert \eta \vert )^{ m-\vert \alpha\vert }
(1 + \vert \sigma \vert )^{\mu-\vert \gamma\vert },\nonumber\\
(s,x,y)\in I^p\times I^{n}, (\xi , \eta, \sigma )\in {\Bbb R}^{n}\times
{\Bbb R}^p.
\end{eqnarray*}
\end{definition}

\noindent A symbol $a\in S^{m,\mu}(I^{n} \times {\Bbb R}^{n}, {\Bbb R}^{p},
{\cal L}({\Bbb C}^r))$ defines an operator
$A$ from $C^{\infty}_c(I^n, {\Bbb C}^r)$ to
$C^{\infty}(I^n, {\Bbb C}^r)$ by the formula
\begin{eqnarray}
\label{loc200}
Au(x,y) & = & (2\pi)^{-2p-q} \int  e^{i[(x-x'-s)\xi+(y-y')\eta+s\sigma]}
a(s,x,y,\xi,\eta,\sigma)\nonumber\\ 
& & u(x',y')\,ds\,dx'\,dy'\,d\xi\,d\eta\,d\sigma,
\end{eqnarray}
\noindent where \(u \in C^{\infty}_{c}(I^{n}, {\Bbb C}^r), x \in I^{p},
y \in I^{q}\).

If $\kappa: I^p\times I^q\rightarrow U=
\kappa(I^p\times I^q)\subset M,
\kappa': I^p\times I^q\rightarrow U'=\kappa'(I^p\times I^q)
\subset M$, are two compartible foliated charts
on $M$, then we can transfer an operator $A$ of
the form (\ref{loc200}) to an operator
$A': C^{\infty}_c(U,E)\rightarrow C^{\infty}(U',E)$.
If, in addition, the kernel of the operator $A'$ is compactly
supported in $U\times U'$, then the operator $A'$ maps
$C^{\infty}_c(U,E)$ to $C^{\infty}_c(U',E)$, and we can
prolong it in a trivial way to an operator
$A':C^{\infty}(M,E)\rightarrow C^{\infty}(M,E)$.
We say that the operator $A'$ obtained in such a way is
{\bf an elementary operator} of class $\Psi ^{m,\mu}(M,{\cal F},E)$.

\begin {definition} The  class $\Psi ^{m,\mu}(M,{\cal F},E)$ 
consists of operators $A$, acting from $C^{\infty}(M,E)$ to 
$C^{\infty}(M,E)$, such that
$A$ can be represented in the form $A=\sum_{i=1}^k A_i + K$,
where $A_i$ are elementary operators
 of class $\Psi ^{m,\mu}(M,{\cal F},E)$, corresponding
to some pairs $\kappa_i,\kappa'_i$ of compartible
foliated charts, $K\in \Psi ^{-\infty}(M,E)$.
\end{definition}

\begin{remark} In notation of \cite{Gu-U,Gr-U}, the Schwartz
kernel $K_A$ of an operator $A\in\Psi ^{m,\mu}(M,{\cal F},E)$ belongs
to the class $I^{m-p/2,\mu+p/2}(M\times M, \Delta, G_{{\cal F}_N})$.
\end{remark}

\begin{example}  
Any tangential pseudodifferential operator
$B\in \Psi^{\mu}({\cal F},E)$ \cite{tang} belongs to the class
$\Psi^{0,\mu}(M,{\cal F},
E)$. Moreover, if $B$ is given by the complete symbol $b(x,y,\xi)$ in some
foliated coordinate chart $\kappa$, then, in the
foliated chart $W(\kappa,\kappa)$,
 the Schwarz kernel of $B$ is represented in the form (\ref{loc200})
with
$$
a(s,x,y,\xi,\eta,\sigma)=b(x,y,\sigma) 
$$
(for a global description, see also Example~\ref{next}).
\end{example}
\begin{example}
\label{next}
Any operator
$C\in \Psi^{m}(M,E)$ belongs to $\Psi^{m,0}(M,{\cal F},E)$.
Moreover, if $C$ is given by the complete symbol $c(x,y,\xi,\eta)$ in some
foliated coordinate chart $\kappa$, then, in the foliated
chart $W(\kappa,\kappa)$, the Schwarz kernel of $C$ is represented in the form
(\ref{loc200}) with
$$
a(s,x,y,\xi,\eta,\sigma)=c(x,y,\xi,\eta).
$$
\end{example}
\begin{example} 
It can be easily seen that two definitions of 
classes $\Psi^{m,-\infty}(M,{\cal F},E)$ are equivalent,
that is,
$$
\Psi^{m,-\infty}(M,{\cal F},E)=\bigcap_{\mu}\Psi^{m,\mu}(M,{\cal F},E).
$$
\end{example}

The operators of class $\Psi^{m,\mu}(M,{\cal F},E)$ can be considered as a
usual pseudifferential operators of order $m+\mu$
with the complete symbol, singular on
the punctured conormal bundle $\tilde{N}^*{\cal F}$. Let us briefly
mention about the corresponding symbolic calculus referring
to \cite{Gu-U,Gr-U} for details.

Let $A$ be an elementary operator of class $\Psi ^{m,\mu}(M,{\cal F},E)$,
given by the formula (\ref{loc200}). Then the principal symbol of $A$ is
a function $\sigma_0(A)$ on $T^*M\backslash N^*{\cal F}$, given locally
by the formula
$$
\sigma_0(A)(x,y,\xi,\eta)= a_{m,\mu}(0,x,y,\xi,\eta,\xi), \xi\not=0,
$$
where $a_{m,\mu}$ is the bihomogeneous component of the complete
symbol $a$ of degree $m$ in $(\xi,\eta)$ and of degree $\mu$ in $\sigma$.

We say that an operator  $A\in\Psi ^{m,\mu}(M,{\cal F},E)$ is elliptic,
if $\sigma_0(A)$ is invertible on $T^*M\backslash N^*{\cal F}$. By
\cite[Proposition 6.4]{Gu-U}, any elliptic operator 
$A\in\Psi ^{m,\mu}(M,{\cal F},E)$ has a parametrix, i.e. an operator
$P\in\Psi ^{-m,-\mu}(M,{\cal F},E)$ such that
\begin{equation}
\label{parametr}
AP=I-R_1, PA=I-R_2,
\end{equation}
where $R_j\in \Psi ^{-1,0}(M,{\cal F},E)+\Psi ^{0,-\infty}(M,{\cal F},E)$,
$j=1,2$.

\begin{proposition}
\label{13}
(1)\ If $A\in  \Psi ^{m_{1},\mu_1}(M,{\cal F},E)$
and $B\in  \Psi ^{m_{2},\mu_2}(M,{\cal F},E)$,
then \(C = AB \in  \Psi ^{m_{1}+m_{2},\mu_1+\mu_2}(M,{\cal F},E)\) and
$\sigma_0(C)=\sigma_0(A)\sigma_0(B)$.

(2)\ Any operator $A\in  \Psi ^{m,\mu}(M,{\cal F},E)$
defines a continuous mapping
$$
A: H^{s,k}(M,{\cal F},E)\rightarrow H^{s-m,k-\mu}
(M,{\cal F}, E)
$$
for any $s$ and $k$.
\end{proposition}

\begin{proof} 
1) A proof is given in \cite{A-U} (see
also \cite[Proposition 1.39]{Gr-U}).

2) As usual, it suffices to consider the case when $A$ is an elementary
operator. Using the standard description of
the Sobolev space  $H^{s,k}(M,{\cal F},E)$ by means of elliptic operators
of class	$\Psi^{s,k}(M,{\cal F},E)$ (we may use the local description,
using the operator $(1+D^2_x+D^2_y)^{s/2}(1+D^2_x)^{k/2}$)
and a parametrix for elliptic operators (\ref{parametr}), 
we can reduce the problem to the
case $s=k=0$ and $\max(m,m+\mu)\leq 0$.
$L^2$ boundedness of operators from  $\Psi ^{m,\mu}(M,{\cal F},E)$ with
 $\max(m,m+\mu)\leq 0$ can be stated by imitating of the H\"ormander's
proof of $L^2$ boundedness of zero-order pseudodifferential operators.
For details see \cite[Theorem 3.3]{Gr-U}.
\end{proof}

\subsection{Symbolic properties of $\Psi ^{m,-\infty}(M,{\cal F},E)$}
Here we turn to more elaborate symbolic properties of classes
$\Psi ^{m,-\infty}(M,{\cal F},E)$. While the algebraic symbolic 
properties of the section follow directly from
the corresponding properties of the algebras ${\cal R}_{\Sigma}$ of
\cite{GS79} (see \cite[Section 3]{GS79}), the presence
of the Sobolev space scale in our geometric case provides norm
estimates in addition to the algebraic symbolic results of \cite{GS79}.

Recall that the principal symbol $p_m$ of 
an operator $P\in \Psi^m(M,E)$ is a smooth section of the vector bundle
${\cal L}(\pi^*E)$ on $\tilde{T}^*M$, where
$\pi: T^*M\rightarrow M$ is the natural projection.

\begin{definition} {\bf The  transversal  principal  symbol}
$\sigma _{P}$  of  an operator $P\in  \Psi ^{m}(M,E)$ is the restriction
of its principal symbol $p_m$ on $\tilde{N}^{\ast}{\cal F}$.
\end{definition}

\begin{proposition}
\label{module}
If $A\in  \Psi ^{m_{1}}(M,E)$
and $B\in  \Psi ^{m_{2},-\infty}(M,{\cal F},E)$,
then $AB$ and $BA$ in  $\Psi ^{m_{1}+m_{2},-\infty}(M,{\cal F},E)\) and
\begin{eqnarray*}
\sigma_{AB}(\gamma,\eta)&=&\sigma_A(\eta)\sigma_B(\gamma,\eta),
(\gamma,\eta)\in G_{{\cal F}_N},\\[6pt]
\sigma_{BA}(\gamma,\eta)&=&
\sigma_B(\gamma,\eta)\sigma_A(dh_{\gamma}^*(\eta)),
(\gamma,\eta)\in G_{{\cal F}_N}.
\end{eqnarray*}
\end{proposition}
\begin{proof} This Proposition follows from the composition
theorem of Fourier integral operators (see, for instance, \cite{H4}).
\end{proof}
From now on, we will assume that $E$ is {\bf holonomy equivariant}, that is,
there is an isometrical action
$$
T(\gamma):E_x\rightarrow E_y, \gamma \in G, \gamma:x\rightarrow y
$$
of the holonomy groupoid $G$ in fibres of $E$.
We have an inclusion $C^{\infty}_c(G)\subset
C^{\infty}_c(G,{\cal L}(E))$, given by $k(\gamma)\mapsto k(\gamma)T(\gamma)$,
and, by (\ref{tang}),
a representation $R_E$ of the algebra $C^{\infty}_c(G)$ in $L^2(M,E)$,
which is a $*$-representation, if the density $dx$ is holonomy invariant.
For any $k\in C^{\infty}_c(G)$, the operator $R_E(k)$
is given by the formula
\begin{equation}
\label{K}
R_E(k)u(x) =\int_{ G^{x}}  k(\gamma ) T(\gamma)u(s(\gamma )) d\lambda ^{x}
(\gamma ), x\in  M, u\in C^{\infty}(M,E).
\end{equation}
Our norm estimates will be given in terms of the following seminorms on
$C^{\infty}_c(G)$:
\begin{equation}
\label{semi}
\|k\|_{s,t,t-l}= \| R_E(k):H^{s,t}(M,{\cal F},E)
\rightarrow H^{s,t-l}(M,{\cal F},E)\|, 
\end{equation}
where $k\in C^{\infty}_c(G)$ and
$s\in {\Bbb R}, t\in {\Bbb R}, l\in {\Bbb R}$. There are defined
the corresponding functional classes:

\begin{definition} The class ${\rm OP}^l(G)$ consists
of all distributions $k\in {\cal D}'(G)$ such that
the operator $R_E(k)$ defines a continuous mapping
$$
R_E(k):H^{s,t}(M,{\cal F},E)\rightarrow H^{s,t-l}(M,{\cal F},E)
$$
for any real $s$ and $t$.
\end{definition}

Recall that an open subset $U$ in $T^*M$ is {\bf a conic neighborhood}
of $N^*{\cal F}$, if
$U$ is a neighborhood of $N^*{\cal F}$, which is invariant under 
the action of ${\Bbb R}_+$ by multiplication.
It is clear that a basis of conic neighborhoods of $N^*{\cal F}$ is formed
by sets $\kappa(U_{\varepsilon}),\varepsilon>0$,
where $U_{\varepsilon}$ is given by
$$
U_{\varepsilon}=\{ (x,y,\xi,\eta)\in I^n\times {\Bbb R}^n:
|\xi|<\varepsilon |\eta|\},
$$
and $\kappa: I^p\times I^q\rightarrow M$ is a foliated coordinate chart.

\begin{definition} 
We say that an operator
$P\in \Psi^l(M,E)$ has {\bf transversal order} 
$m\leq l$ ($P\in \Psi^m(N^*{\cal F},E)$),
if $P$ has order $m$ in some conic neighborhood of $N^*{\cal F}$,
that is, its complete symbol $p$ in any
foliated coordinate system satisfies the following condition:
there is a $\varepsilon>0$
such that, for any multiindices $\alpha$, $\beta$,
there is a constant $C_{\alpha,\beta }>0$ such that
\begin{equation}
\label{tra}
\vert \partial^{\alpha}_{(\xi,\eta)} \partial^{\beta}_{(x,y)}
p(x,y,\xi,\eta )\vert
\leq C_{\alpha \beta}(1 +\vert \xi \vert + \vert \eta \vert )^
{ m-\vert \alpha \vert }, 
(x,y,\xi,\eta)\in U_{\varepsilon}.
\end{equation}
\end{definition}

The main fact, which relates the notion of transversal order with
the classes $\Psi^{m,\mu}(M,{\cal F},E)$, consists in the following
inclusion:
\begin{equation}
\label{incl}
\Psi^l(M,E)\bigcap \Psi^m(N^*{\cal F},E)\subset
\Psi^{m,l-m}(M,{\cal F},E).
\end{equation}
Indeed, it can be easily checked by a straightforward calculation, that,
if $p\in S^l(I^n\times {\Bbb R}^n)$ satisfies (\ref{tra}), then
$p(x,y,\xi,\eta)(1+|\xi|^2)^{m-l}$ belongs to
$S^m(I^n\times {\Bbb R}^n)$, from where (\ref{incl}) follows immediately
(see also \cite[Theorem 18.1.35]{H3}).

By Proposition~\ref{13} and  (\ref{incl}), 
any operator $P\in\Psi^m(N^*{\cal F},E)\bigcap \Psi^l(M,E)$
defines a continuous mapping
\begin{equation}
\label{cont}
P:H^{s,k}(M,{\cal F},E)\rightarrow H^{s-m,k-l+m}(M,{\cal F},E),
\end{equation}
that, in its turn, gives immediately the following proposition.

\begin{proposition}
\label{tr}
Suppose that an operator $P\in \Psi^l(M,E)$ has transversal order $m\leq l$.
 Then, for any $k\in
C^{\infty}_c(G,{\cal L}(E))$, the operators $R_E(k)P$ and $PR_E(k)$ belong
to $\Psi ^{m,-\infty}(M,{\cal F},E)$ with the following norm
estimates
\begin{eqnarray*}
\|R_E(k)P:H^{s,t}(M,{\cal F},E)\rightarrow H^{s-m,r}(M,{\cal F},E)\|
&\leq &C\|k\|_{s-m,t-l+m,r},\\
\|PR_E(k):H^{s,t}(M,{\cal F},E)\rightarrow H^{s-m,r}(M,{\cal F},E)\|
&\leq &C\|k\|_{s,t,r+l-m}.
\end{eqnarray*}
\end{proposition}

We will denote by $\hbox{ad}\, T(\gamma)$  the action of $G$
in fibres of the bundle ${\cal L}(\pi^*E)$, induced by $T(\gamma)$.

\begin{definition}
We say that the transversal principal symbol  of
an operator $P\in\Psi^m(M,E)$ is {\bf holonomy invariant}, if,
for any  smooth  leafwise path $\gamma $ from $x$ to $y$,
the following equality holds:
$$
{\rm ad}\,T(\gamma)[\sigma _{P}(dh_{\gamma}^{\ast}(\xi ))] =
\sigma _{P}(\xi ), \xi  \in N^{ \ast}_{y}{\cal F}.
$$
\end{definition}

\begin{remark} 
The assumption of existence of positive order
pseudodifferential operators with the
holonomy invariant transversal principal symbol implies rather strong
restrictions on the foliated manifold under consideration. There are
examples of such operators on every Riemannian foliation given
by the transverse signature operator, and, in fact, this assumption is
equivalent to a slightly more general assumption on the foliation to be
transversally Finsler (as introduced in \cite{tang}). For general foliations,
one can use a generalized notion of holonomy invariance, based on
more sophisticated transversal pseudodifferential calculus. As an example, we
point out the treatment of triangular
Riemannian manifolds, based on hypoelliptic operators and Beals-Greiner
pseudodifferential calculus \cite{Co-M}.
We will discuss this subject in more details in \cite{noncom2}.
\end{remark}

By Proposition~\ref{module}, if $P$ is a pseudodifferential operator from
$\Psi ^{m}(M,E)$ with the holonomy invariant transversal principal symbol, and
$k\in C^{\infty}_c(G)$, the operator  $[P,R_E(k)]$ belongs to the class
$\Psi^{m-1,-\infty}(M,{\cal F},E)$.
Moreover, we have the following norm estimate:
\begin{proposition}
\label{com}
Let $P$ be a pseudodifferential operator from $\Psi ^{m}(M,E)$ with the
holonomy invariant transversal principal symbol. Then, for any
$k\in {\rm OP}^l(G)$,
the operator $[P,R_E(k)]$ defines a continuous map
$$
[P,R_E(k)]:H^{s,t}(M,{\cal F},E)\rightarrow H^{s-m+1,t-l-1}(M,{\cal F},E)
$$
with the following norm estimate
\begin{eqnarray*}
\|[P,R_E(k)] & : & H^{s,t}(M,{\cal F},E)\rightarrow H^{s-m+1,t-l-1}
(M,{\cal F},E)\|\nonumber\\
& & \leq C\max (\|k\|_{s,t,t-l}, \|k\|_{s-m+1,t-1,t-l-1}).
\end{eqnarray*}
\end{proposition}
\begin{proof} Let $p\in S^m(I^n,{\cal L}({\Bbb C}^r))$ be
the complete symbol of the operator $P$ in some foliated
chart with a complete asymptotic expansion
$p\sim \sum_{j=0}^{\infty} p_{m-j}$, 
$p_{m-j}(x,y,\xi,\eta)$ is homogeneous in $(\xi,\eta)$
of degree $m-j$. By the holonomy invariance assumption, we can choose a
trivialization of the bundle $E$ so that the transversal
principal symbol in this coordinate system will be
a matrix-valued function, independent of $x$:
$$
p_m(x,y,0,\eta)=p_m(y,\eta).
$$
Using the Taylor formula, we represent $p_m$ in the form
$$
p_m(x,y,\xi,\eta)=p_m(y,\eta)+\sum_{i=1}^p p_{m,i}(x,y,\xi,\eta)\xi_i,
$$
where $p_{m,i}$ are homogeneous in $(\xi,\eta)$
of degree $m-1$.

Let $P_1$ be an operator with the complete symbol
$p_m(y,\eta)$, $R_1$ be an operator with the complete symbol
$\sum_{i=1}^p p_{m,i}(x,y,\xi,\eta)\xi_i$. Gluing together
these local operators into global ones in a standard way,
we get a representation
$P=P_1+R_1$,
where $P_1\in \Psi ^{m}(M,E)$, $[P_1,R_E(k)]\in
\Psi ^{-\infty}(M,E)$, $R_1\in \Psi ^{m-1,1}(M,{\cal F},E)$.
So we have $[P,R_E(k)]=[R_1,R_E(k)]\hbox{\rm mod}\ \Psi ^{-\infty}(M,E)$,
from where Proposition~\theproposition\ is immediate.
\end{proof}

\subsection{Residue trace}
Now we turn to a trace extension	and a Wodzicki type residue
for operators
of class $\Psi ^{m,-\infty}(M,{\cal F},E)$. These results are particular
cases of results \cite{Gu_1,Gu_2} on Fourier integral operators,
but we make use of the local structure of operators in 
question given by (\ref{loc1}) to derive 
directly all necessary facts

For any $\sigma\in C^{\infty}({\Bbb R}^q\backslash \{0\})$, 
homogeneous of order
$q$, i.e.
$\sigma(\lambda\eta)=\lambda^d\sigma(\eta)$
for any $\eta\not= 0$ and $\lambda \in {\Bbb R}^*_+$,
let
$$
S({\sigma})=\int_{|\eta|=1} \Tr \sigma(\eta)d\eta.
$$

For any function $\phi$ on ${\Bbb R}^q\backslash \{0\}$, let
$$
\phi_{\lambda}(\eta)=\lambda^q\phi(\lambda\eta),
\lambda>0, \eta\in {\Bbb R}^q\backslash \{0\}.
$$

Recall the following fact on continuation of a homogeneous smooth
function on ${\Bbb R}^q\backslash \{0\}$ to a homogeneous distribution
in ${\Bbb R}^q$, see \cite{Ho}, Theorems 3.2.3 and 3.2.4.
\begin{lemma}
\label{hom}
Let $\sigma\in C^{\infty}({\Bbb R}^q\backslash \{0\})$ be 
homogeneous of order $d$ in $\eta \in {\Bbb R}^q$.
\medskip
\par
\noindent (1)\ If $d\not \in \{-q-k:k\in {\Bbb N}\}$,
$\sigma$ extends to a homogeneous distribution $\tau$
on ${\Bbb R}^q$.
\medskip
\par
\noindent (2)\ If $d=-q-k$, there is an extension $\tau$
of $\sigma$, satisfying the condition
$$
(\tau,\phi)=\lambda^{-q-k}(\tau,\phi_{\lambda})+
\log \lambda \sum_{|\alpha|=k} S(\eta^{\alpha}\sigma)
\partial_{\eta}^{\alpha}\phi(0)/\alpha!, \alpha>0.
$$
In particular, the obstruction to an extension of $\sigma\in
{\cal D}'({\Bbb R}^q)$, homogeneous in $\eta$, is given by
$S(\eta^{\alpha}\sigma)$, $|\alpha|=k$.
\end{lemma}

Let a functional $L$ be given by the formula
$$
L(\sigma)=(2\pi)^{-q}\int \Tr \sigma(\eta)d\eta,
$$
which is well-defined on symbols $\sigma\in S^m_{cl}({\Bbb
R}^q)$
of order $m<-q$.
\begin{lemma}[\cite{KV1,Co-M}]
The functional $L$ has an unique holomorphic extension
$\tilde{L}$ to the space of classical symbols
$S^z_{cl}({\Bbb R}^q)$ of non-integral order $z$.
The value of $\tilde{L}$ on a symbol
$\sigma\sim \sum \sigma_{z-j}$ is given by
$$
\tilde{L}(\sigma)=
(2\pi)^{-q}\int  \Tr (\sigma -\sum_0^N \tau_{z-j}) d\eta,
$$
where $\tau_{z-j}$ is the unique homogeneous extension
of $\sigma_{z-j}$, given by Lemma~\ref{hom}, $N\geq {\rm Re}\,z
+q$.
\end{lemma}

The trace of a pseudodifferential operator
$A\in\Psi^{m,-\infty}(I^n,I^p,{\Bbb C}^r)$ given by (\ref{loc})
with $m<-q$ is given by the formula
$$
\tr(A)=(2\pi)^{-q}\int  \Tr k(x,x,y,\eta) dx dy d\eta.
$$
The following formula provides an extension of the trace
to a pseudodifferential operator $A\in
\Psi^{z,-\infty}(I^n,I^p,{\Bbb C}^r)$ of arbitrary non-integral order
$z\in {\Bbb C}\backslash {\Bbb Z}$:
$$
\TR(A)=(2\pi)^{-q}\int  \Tr \tilde{L}(k(x, x ,y,\eta))dx dy d\eta.
$$
This definition can be extended to elementary operators, and, by 
linearity, to all operators
$P\in \Psi^{z,-\infty}(M,{\cal F},E), z\in {\Bbb C}\backslash {\Bbb Z}$.

If an operator $A\in\Psi^{m,-\infty}(I^n,I^p,{\Bbb C}^r)$ is written
in the form (\ref{loc1}), then $\TR(A)$ is given by
\begin{equation}
\label{loc100}
\TR(A)=\int  \TR(A(x,x))\,dx,
\end{equation}
where $\TR(A(x,x))$ denotes the extension of the usual trace of the
pseudodifferential operator $A(x,x)$ on $I^q$, defined in \cite{KV1,KV2}.
Using (\ref{loc100}) and \cite{KV1,KV2}, we immediately obtain
the following proposition.

\begin{proposition}
The linear functional $\TR$ on the class $\Psi^{\alpha,-\infty}(M,{\cal F},E)$
of classical pseudodifferential operators of
orders $\alpha\in m+{\Bbb Z},m\in {\Bbb C}\backslash {\Bbb Z}$,
has the following properties:
\medskip
\par
\noindent (1) It coincides with the usual trace $\tr$
for ${\rm Re}\, \alpha<-q$.
\medskip
\par
\noindent (2) It is a trace functional, i.e.
$\TR([A,B])=0$ for any $A\in \Psi^{\alpha_1,-\infty}_{cl}
(M,{\cal F},E)$ and $B\in \Psi^{\alpha_2,-\infty}_{cl}
(M,{\cal F},E)$, $\alpha_1+\alpha_2\in m+{\Bbb Z}$.
\end{proposition}

Now let us turn to the residue trace of operators from
$\Psi^{m,-\infty}(M,{\cal F},E)$. As above, it suffices to define
the residue trace for an elementary operator.
Given an operator $A\in\Psi^{m,-\infty}(I^n,I^p,{\Bbb C}^r)$,
we define its residue form $\rho_A$ as
$$
\rho_A = \Tr k_{-q}(x,x,y,\eta) dx dy d\eta,
$$
and the residue trace $\tau(A)$ as
$$
\tau(A)=\int_{|\eta|=1}\Tr k_{-q}(x,x,y,\eta) dx dy d\eta.
$$
If an operator $A\in\Psi^{m,-\infty}(I^n,I^p,{\Bbb C}^r)$
is written in the form (\ref{loc1}),
then its residue trace $\tau(A)$ is given by
\begin{equation}
\label{residue}
\tau(A)=\int  \tau(A(x,x))\,dx,
\end{equation}
where $\tau(A(x,x))$ denotes the residue trace of the
pseudodifferential operator $A(x,x)$ on $I^q$ due to \cite{Gu85,Wo}.

Using (\ref{residue}), it can be easily checked that, for any
$A\in \Psi^{m,-\infty}(M,{\cal F},E)$,
its residue form $\rho_A$ is an invariantly defined form on
$\tilde{N}^*{\cal F}$, and the residue
trace $\tau(A)$ is given by integration of the residue form $\rho_A$
over the spherical conormal bundle $SN^*{\cal F}=\{\nu \in N^*{\cal F}
: |\nu|=1\}$:
$$
\tau(A)=\int_{SN^*{\cal F}}\rho_A.
$$

\begin{remark}
Let $A\in \Psi^{m,-\infty}(M,{\cal F},E)$, and $K_A\in I^{m-p/2}
(M\times M, G'_{{\cal F}_N})$ be its Schwarz kernel. In terms
of \cite{Gu_1}, the residue trace
$\tau(A)$ of $A$ is defined as the residue pairing of $K_A$ with the
delta function $\delta_{\Delta}$ of the diagonal $\Delta$ in
$T^*M$. The holonomy groupoid $G_{{\cal F}_N}$ and
$\Delta$ are intersected cleanly in the $p+2q$-submanifold
$G^{(0)}_{{\cal F}_N}=N^*{\cal F}$, and, therefore, this residue
pairing is defined by \cite[Theorem 4.1]{Gu_1}.
The main properties, concerning to the trace extension and the residue form
for general Fourier integral operators, are due to \cite[Theorem 2.1]{Gu_1}.
\end{remark}

Now we relate the trace extension and the residue trace for operators
of $\Psi^{m,-\infty}(M,{\cal F},E)$.
At first, let us give a definition of a holomorphic family of
pseudodifferential operators of class $\Psi^{*,-\infty}(M,{\cal F},E).$
As usual, it is sufficient to do this for elementary operators.

\begin{definition} We say that a family $A(z)\in 
\Psi^{f(z),-\infty}_{cl}(I^{n}, I^{p}, {\Bbb C}^r)$ is holomorphic
(in a domain $D\subset {\Bbb C}$), if:
\medskip
\par
\noindent (1) the order $f(z)$ is a holomorphic function;
\medskip
\par
\noindent (2) $A(z)$ is given by a classical symbol $k(z)\in S^{f(z),-\infty}
(I^{p}\times I^{p}\times I^{q} \times {\Bbb R}^{q},
{\cal L}({\Bbb C}^r))$, represented as an asymptotic sum
$$
k(z,x,x',y,\eta)\sim \sum_{j=0}^{\infty} \theta(\eta)
k_{f(z)-j}(z,x,x',y,\eta),
$$
which is uniform in $z$, and the homogeneous components
$k_{z-j}(z,x,x',y,\eta)$ are holomorphic in $z$.
\end{definition}

\begin{proposition}
\label{extension}
For any holomorphic family
$A(z)\in \Psi^{m + {\rm Re}\,z,-\infty}_{cl}(M,{\cal F},E), z\in D\subset
{\Bbb C}$, the function $z\mapsto \TR(A(z))$ is meromorphic with no
more than simple poles at $z_k=-m-q+k\in D\bigcap {\Bbb Z},k\geq 0$ and
with
$$
\res_{z=z_k} \TR(A(z))=\tau(A(z_k)).
$$
\end{proposition}
\begin{proof} The proposition is an immediate consequence of
(\ref{residue}) and of the similar fact for usual pseudodifferential
operators \cite{KV1,KV2}.
\end{proof}

\section{Transversally elliptic operators}
\subsection{Definition and basic properties}
As above, we assume that $M$ is a closed foliated manifold, and $E$ is
a holonomy equivariant Hermitian vector bundle $E$ on $M$.

\begin{definition} 
We say that an operator
$P\in \Psi^m(M,E)$ is {\bf transversally elliptic},
if the transversal principal symbol of $P$
is invertible for any $\xi\in \tilde{N}^*{\cal F}$.
\end{definition}

\begin{remark}
\label{transell}
The condition of transversal
ellipticity implies that the transversal principal symbol of $P$
is invertible for any $\xi$
in some conic neighborhood of $N^*{\cal F}$.
It also implies that, in any foliated coordinate system,
there are $\varepsilon>0$ and $c>0$ such that
$$
|(p_m(x,y,\xi,\eta )v,v)|\geq
c(1 +\vert \xi \vert + \vert \eta \vert )^{m}\|v\|,
(x,y,\xi,\eta)\in U_{\varepsilon}, v\in {\Bbb C}^r,
$$
where $p_m$ is the homogeneous component of degree $m$ in $(\xi,\eta)$
of the complete symbol of $P$.
\end{remark}

Performing the standard parametrix construction in the conic neighborhood
$U_{\varepsilon}$ of $N^*{\cal F}$ given by Remark~\ref{transell}, it is easy
to get the following proposition (see also \cite{trans}).
\begin{proposition}
For a transversally elliptic  operator $P\in
\Psi ^{m}(M,E)$, there  exists  a
parametrix, that is, an operator $Q\in  \Psi ^{-m}(M,E)$
such that
\begin{equation}
\label{param}
PQ = I - R_{1}, QP = I - R_{2},
\end{equation}
\noindent where $R_{j} \in
\Psi ^{0}(M,E)\bigcap \Psi^{-\infty}(N^*{\cal F},E), j = 1,2$.
\end{proposition}

Existence of a parametrix implies a transverse elliptic regularity theorem
in a usual manner:
\begin{proposition}
\label{reg}
Given a transversally elliptic operator $P\in  \Psi ^{m}(M,E)$, and a section
$u$ such that
$u \in H^{s+m-N,k+N}(M,{\cal F},E)$ for some $N>0$
and $Pu \in H^{s,k}(M,{\cal F},E)$, we have $u\in H^{s+m,k}(M,{\cal F},E)$
and
$$
\| u\| _{s+m,k} \leq C(\|Pu\| _{s,k} + \| u\| _{s+m-N,k+N}).
$$
\end{proposition}

\subsection{Complex powers}
\label{complex}
Throughout in this section, we assume that an operator $A\in\Psi^{m}(M,E)$
satisfies the following conditions:
\medskip
\par
(T1) $A$ is a transversally elliptic pseudodifferential operator with
the positive transversal principal symbol;
\medskip
\par
(T2) $A$ is essentially self-adjoint on the initial domain $C^{\infty}(M,E)$,
and its closure is invertible and positive definite as an unbounded operator
in the Hilbert space $L^2(M,E)$.

\begin{remark} The assumption (T2) may be
considered as an equivariance type condition, which is, usually, assumed
for transversally elliptic operators.
\end{remark}

\begin{proposition}
\label{inverse}
Let $A \in \Psi^{m}(M,E)$ be as above. Then, for any
$\lambda\not\in {\Bbb R}_+$, the resolvent operator
$(A-\lambda )^{-1}$ is  represented as
\begin{equation}
\label{rez}
(A-\lambda )^{-1}= P(\lambda) + R_1(\lambda )(A-\lambda )^{-1},
(A-\lambda )^{-1}= P(\lambda) + (A-\lambda )^{-1}R_2(\lambda ),
\end{equation}
where:

(1)\ $P(\lambda ) \in \Psi^{-m}(M,E)$ is an operator,
which complete symbol in any foliated coordinate
system is supported in some conical neighborhood $U_{\varepsilon}$
of $N^*{\cal F}$ and satisfies the estimates
\begin{eqnarray*}
\vert D^{\beta}_{(x,y)}D^{\alpha}_{(\xi,\eta)}p(x,y,\xi ,\eta ,\lambda )\vert
&\leq &C_{\alpha \beta }(1+\vert \xi \vert +\vert\eta\vert+\vert 
\lambda \vert ^{1/m})^{-m}\\[6pt]
& &(1+\vert \xi \vert +\vert \eta \vert )^{- \vert \alpha\vert},
(x,y,\xi,\eta)\in U_{\varepsilon}, \lambda  \in  \Lambda_{\delta}
\end{eqnarray*}
for any $\delta  > 0$ and for any multi-indices $\alpha$ and $\beta$;

(2)\ $R_j(\lambda ) \in \Psi^{0}(M,E)\bigcap \Psi^{-\infty}(N^*{\cal F},E),
j=1,2,$
with the complete symbol $r_j(\lambda)$, satisfying the following estimates:
\begin{eqnarray}
\label{rn}
\vert D^{\beta}_{(x,y)}D^{\alpha}_{(\xi,\eta)}r_j(x,y,\xi ,\eta ,\lambda )\vert
& \leq & C_{\alpha \beta }(1+\vert \xi \vert +|\eta|+
\vert \lambda \vert ^{1/m})^{-m}\nonumber\\[6pt]
& &(1+\vert \xi \vert +\vert \eta \vert )^{- \vert \alpha\vert},
(x,y,\xi,\eta)\in U_{\varepsilon_1}, \lambda  \in  \Lambda_{\delta};
\nonumber\\[6pt]
\vert D^{\beta}_{(x,y)}D^{\alpha}_{(\xi,\eta)}r_j(x,y,\xi ,\eta ,\lambda )\vert
& \leq & C_{\alpha \beta N}(1+\vert \xi \vert +\vert \eta \vert )^
{- \vert \alpha\vert-N+m},\nonumber \\[6pt]
& & (x,y,\xi,\eta)\in I^p\times I^q\times {\Bbb R}^p\times {\Bbb R}^q,
\lambda  \in  \Lambda_{\delta}
\end{eqnarray}
for any $\delta  > 0$, for any natural $N$ and for any multi-indices
$\alpha$ and $\beta$.

Moreover, the principal symbol $p_{-m}(\lambda)$ of $P(\lambda )$ 
is equal to $(a_m-\lambda)^{-1}$ in some conic
neighborhood of $N^*{\cal F}$ with $a_m$, being the principal symbol of $A$.
\end{proposition}

\begin{proof} We will prove the proposition,
performing the standard construction due to Seeley of a
parametrix $P(\lambda)$ for the operator $A-\lambda $ as an operator with
a parameter in some conic neighborhood of $N^*{\cal F}$.

Denote by $\Lambda _{\delta}$ the angle in the complex plane:
$$
\Lambda _{\delta} = \{\lambda  \in  {\Bbb C} : \vert \arg
\lambda \vert  > \delta \}.
$$

Fix some foliated coordinate system.
Let $a \sim \sum^{\infty}_{j=0} a_{j}$ be an asymptotic expansion of the
complete symbol of the operator $A$ in this system. By (T1),
there is a $\varepsilon>0$ such that
$$
a_m(x,y,\xi,\eta)\geq C(1+|\xi|+|\eta|)^m,
(x,y,\xi,\eta)\in U_{\varepsilon}.
$$
For any $\delta > 0$, define
functions $p_{-m-l}(\lambda), \lambda \in \Lambda _{\delta},l=0,1,\ldots,$ in
$U_{\varepsilon}$ by the following system
\begin{eqnarray*}
& (a_{m}-\lambda)p_{-m}=1,& \\
& (a_{m}-\lambda)p_{-m-l}+ \sum_{j < l, j+k+\vert\alpha\vert=l}
{\partial}_{\xi}^{\alpha}
b_{-m-j}D^{\alpha}_{x}a_{m-k}/\alpha! =0,& l>0.
\end{eqnarray*}

\noindent It can be easily
checked that the functions $p_{-m-l}(\lambda)$ satisfy the following estimates
\begin{eqnarray*}
\vert D^{\beta}_{(x,y)}D^{\alpha}_{(\xi,\eta)}
p_{-m-l}(x,y,\xi ,\eta ,\lambda )\vert & \leq &
C_{\alpha \beta }(1+\vert \xi \vert +\vert\eta 
\vert+\vert \lambda \vert ^{1/m})^{-m}\\[6pt]
& &(1+\vert \xi \vert +\vert \eta \vert )^{- \vert \alpha\vert}, 
(x,y,\xi,\eta)\in U_{\varepsilon}, \lambda  \in  \Lambda_{\delta},
\end{eqnarray*}
\noindent where $\alpha$ and $\beta$ are any multi-indices.
Take $p$ as an asymptotic sum of symbols with a parameter:
\(p \sim \sum_{j=0}^{+\infty} p_{-m-j}\).
Then $p$ satisfies the estimates
\begin{eqnarray}
\label{pn}
\vert D^{\beta}_{(x,y)}D^{\alpha}_{(\xi,\eta)}p(x,y,\xi ,\eta ,\lambda )\vert &
\leq &C(1+\vert \xi \vert +\vert\eta\vert+\vert \lambda \vert ^{1/m})^{-m}
(1+\vert \xi \vert +\vert \eta \vert )^{- \vert \alpha\vert}, \nonumber \\[6pt]
& &(x,y,\xi,\eta)\in U_{\varepsilon}, \lambda  \in  \Lambda_{\delta}.
\end{eqnarray}

Let $\theta_0\in C^{\infty}_c({\Bbb R})$
such that ${\rm supp}\, \theta_0\subset 
(-\varepsilon, \varepsilon)$, $\theta_0(\tau)=1$ for any $\tau\in
(-\varepsilon_1, \varepsilon_1)$ with some
$\varepsilon_1<\varepsilon$, and $\theta\in C^{\infty}
({\Bbb R}^p\times{\Bbb R}^q)$ be given by
$\theta(\xi,\eta)=\theta_0(\eta/\xi)$, if $|\eta|<\varepsilon|\xi|$,
and $\theta(\xi,\eta)=0$ in the opposite case.
Let us take a covering of $M$ by foliation charts, construct in any foliation
patch of this covering an operator with the 
complete symbol $\theta p(\lambda)$,
and glue these local operators in a global operator $P(\lambda ) \in
\Psi^{-m}(M,E), \lambda  \in  \Lambda_{\delta}$ by means of a partition
of unity. It can be easily seen that
\begin{equation}
\label{pd}
P(\lambda )(A-\lambda ) = I - R_1(\lambda ),
(A-\lambda )P(\lambda) = I - R_2(\lambda ), \lambda  \in  \Lambda_{\delta},
\end{equation}
\noindent where $R_j(\lambda ) \in \Psi^{0}(M,E)\bigcap
\Psi^{-\infty}(N^*{\cal F},E), j=1,2,$ has the complete
symbol $r_j(\lambda)$, satisfying the estimates (\ref{rn}).
By (T2), the operator  $A-\lambda $  is
invertible as  an  unbounded  operator in $L^2(M,E)$
for all $\lambda \not\in {\Bbb R}_+$ with
the following estimate for the norm of its inverse:
$$
\|(A-\lambda)^{-1}\|\leq C/|\lambda|.
$$
Using (\ref{pd}), we get the representation~(\ref{rez})
for the resolvent, that completes the proof.
\end{proof}

\begin{proposition}
\label{norm}
Let $A \in \Psi^{m}(M,E)$ be as above. Then, for any
$\lambda\not\in {\Bbb R}_+$, the resolvent operator
$(A-\lambda )^{-1}$ can be  represented as
\begin{equation}
\label{rez1}
(A-\lambda )^{-1}= P(\lambda) + T(\lambda),
\end{equation}
\noindent where:

(1) $P(\lambda ) \in \Psi^{-m}(M,E)$
 satisfies the following norm estimates:
$$
\|P(\lambda ) : H^{s,k}(M,{\cal F},E) \rightarrow  H^{s,k}(M,{\cal F},E)\|
\leq	C_{s,k,\delta}(1+\vert \lambda \vert )^{-1}, \lambda	\in
\Lambda _{\delta},
$$
$$
\|P(\lambda ) : H^{s,k}(M,{\cal F},E) \rightarrow  H^{s+m,k}(M,{\cal F},E)\|
\leq  C_{s,k,\delta}, \lambda  \in  \Lambda _{\delta},
$$
for any $s \in {\Bbb R}$, $k \in {\Bbb R}$ and $\delta>0$;

(2) $T(\lambda)$  satisfies the following norm estimates:
$$
\|T(\lambda ) : H^{t,-t}(M,{\cal F},E) \rightarrow
H^{s,-s}(M,{\cal F},E)\|	\leq	C_{s,t,\delta }(1+\vert \lambda \vert )^{-1},
\lambda  \in  \Lambda _{\delta},
$$
for any $s,t$ and $\delta>0$.
\end{proposition}

\begin{proof} By (\ref{rez}),
$(A-\lambda )^{-1}= P(\lambda) + R_1(\lambda )P(\lambda)+
R_1(\lambda )(A-\lambda )^{-1}R_2(\lambda )$,
so we get (\ref{rez1}) with
\begin{equation}
\label{Tlambda}
T(\lambda)=R_1(\lambda )P(\lambda)+R_1(\lambda )(A-\lambda )^{-1}R_2(\lambda ).
\end{equation}
Since $R_j(\lambda ) \in \Psi^{0}(M,E)\bigcap \Psi^{-\infty}(N^*{\cal F},E),
j=1,2$, by (\ref{cont}), $R_j$ defines a continuous map
from $H^{t,-t}(M,{\cal F},E)$ to $H^{s,-s}(M,{\cal F},E)$ for any $s$ and
$t$, that implies the same is true for $T(\lambda)$.
The desired norm estimates for operators $P(\lambda)$ and
$T(\lambda)$ follow immediately from the symbol estimates (\ref{pn})
and (\ref{rn}).
\end{proof}

Now we turn to a construction of complex powers $A^{z}$ for a transversally
elliptic operator $A \in \Psi^{m}(M,E)$, satisfying to the conditions (T1)
and (T2).

Let $\Gamma$  be a  contour  in the complex plane of the form \(\Gamma =
\Gamma_{1} \bigcup \Gamma_{2} \bigcup \Gamma_{3}\), where
\(\lambda = re^{i\alpha }, +\infty > r > \rho,\) on $\Gamma_{1}$,
\(\lambda = \rho e^{i\phi }, \pi > \phi > -\pi,\) on $\Gamma_{2}$,
\(\lambda = re^{-i\alpha }, \rho < r < +\infty,\) on $\Gamma_{3}$,
\noindent where $\alpha \in (0,\pi)$ is arbitrary,
and the constant $\rho > 0$ is chosen in such a way
that the disk of
the radius $\rho$, centered at the origin, isn't contained in $\sigma(D)$.

A bounded operator $A^{z}$, $\mbox{Re}\, z< 0$, in	$L^{2}(M,E)$	is
defined by the formula

$$
A^{z} = \frac{i}{2\pi} \int_{\Gamma} \lambda ^{z}(A-\lambda )^{-1}d\lambda ,
$$

\noindent where a branch of the analytic function $\lambda ^{z}$ is chosen
so that \(\lambda ^{z} = e^{z\,\ln \lambda}\) for $\lambda > 0$.
This definition is extended to all $z$ by
\begin{equation}
\label{degree}
A^{z} =  A^{z-k}A^{k}
\end{equation}
for any $z$, $\hbox{Re}\,z < k$,
where $k$ is natural and $A^{k}$ is the usual power of the operator $A$.
The following proposition provides a descripiton of the complex powers $A^z$.

\begin{proposition}
\label{powers}
Under the current hypotheses on the operator $A$, the
operator $A^z$ has the form
\begin{equation}
\label{decom}
A^z= P(z)+T(z),
\end{equation}
where $P(z)$ is a holomorphic family of pseudodifferential operators
of class $\Psi^{mz}(M,E)$ with the principal symbol $p_{m{\rm Re}\,z}(z)$, 
being equal to $(a_m)^z$ in some conic neighborhood of $N^*{\cal F}$
($a_m$ is the principal symbol of $A$),
and, for any $s,l$ and $z, {\rm Re}\ z\leq k$,
$T(z)$ defines a continuous mapping
\begin{eqnarray*}
T(z):H^{s,-s}(M,{\cal F},E)&\rightarrow &H^{l,-l-k}(M,{\cal F},E), k>0,\\
T(z):H^{s,-s}(M,{\cal F},E)&\rightarrow &H^{l,-l}(M,{\cal F},E), k\leq 0.
\end{eqnarray*}
\end{proposition}

\begin{proof} A proof of the proposition can be obtained by a
straighforward repetition of the proof of \cite[Proposition 7.3]{trans}.
Namely, let us write $A^z$ for $\mbox{Re}\ z<0$ as $A^z=P(z)+T(z)$, where
$$
P(z) = \frac{i}{2\pi} \int_{\Gamma} \lambda ^{z}P(\lambda )d\lambda ,
T(z) = \frac{i}{2\pi} \int_{\Gamma} \lambda ^{z}T(\lambda )d\lambda ,
$$
$P(\lambda)$ and $T(\lambda)$ are given by Proposition~\ref{norm}. All the
statements about $P(z)$ and $T(z)$ can be easily checked in a standard
way (see \cite{trans} for more details).
\end{proof}

\begin{remark}
\label{adding}
We can get $P(z)$ to be an elliptic operator of class 
$\Psi^{m{\rm Re}\,z}(M,E)$ with the positive principal symbol, adding 
to $P(z)$ an appropriate operator of
class $\Psi^{m{\rm Re}\, z}(M,E)
\bigcap \Psi^{-\infty}(N^*{\cal F},E)$,. 
\end{remark}

\subsection{$G$-trace}
Before going to the distributional zeta-function of transversally
elliptic operators, we introduce a general scheme of defining
distributional spectral invariants for transversally elliptic operators
based on the notion of the $G$-trace and prove some existence results
for such invariants, that may have its own interest.

\begin{definition}
We say that a bounded operator
$T$ in $L^2(M,E)$ is {\bf an operator of $G$-trace class},
if, for any $k\in C^{\infty}_c(G)$, the operator $R_E(k)T$
is a trace class operator, and a functional
${\trG}(T)$ on $C^{\infty}_c(G)$, defined
by the formula
$$
({\trG}(T),k)=\tr R_E(k)T, k\in C^{\infty}_c(G),
$$
is a distribution on $G$.
In this case, the distribution ${\trG}(T)\in {\cal D}'(G)$
is called  {\bf the $G$-trace} of the operator $T$.
\end{definition}

For any integral operator $T$ on $C^{\infty}(M,E)$
with the smooth kernel $K_T$, its $G$-trace ${\trG}(T)$ is
a smooth function on the holonomy groupoid $G$, given
by the formula
$$
{\trG}(T)(\gamma)=K_T(r(\gamma),s(\gamma)),\gamma\in G.
$$
Otherwise speaking, the $G$-trace ${\trG}(T)$ is obtained
by pulling back of the integral kernel $K_T$ via the map
$(r,s):G\rightarrow M\times M$.

By Propositions~\ref{trace} and \ref{13}, we immediately obtain that
any $P\in \Psi^{m,\mu}(M,{\cal F},E)$ with $m<-q$ is an operator
of $G$-trace class.
Indeed, the following, more general proposition is valid.
\begin{proposition}
\label{G-trace}
Let $T$ be a bounded operator in $C^{\infty}(M,E)$, which extends to
a bounded operator from $L^2(M,E)$
to $H^{-m,-\mu}(M,{\cal F},E)$ with some $m<-q$ and
$\mu$. Then $T$ is an operator
of $G$-trace class with
the following estimate of its $G$-trace:
$$
|({\trG}(T),k)|\leq C_{\varepsilon}\|k\|_{q+\varepsilon,-m-\mu-q-\varepsilon,
p+\varepsilon_1},
$$
where $\varepsilon>0$ and $\varepsilon_1>0$ are arbitrary
constants such that $s>q+\varepsilon$.
\end{proposition}
\begin{proof} Let us fix some $\varepsilon>0$
and $\varepsilon_1>0$ such that $m<-q-\varepsilon$. By Proposition~\ref{trace},
we have
$$
|\tr R_E(k)P|\leq C\|R_E(k)P: L^2(M,E)\rightarrow
H^{q+\varepsilon, p+\varepsilon_1}(M,{\cal F},E)\|.
$$
Using the embedding
$H^{-m,-\mu}(M,{\cal F},E)\subset
H^{q+\varepsilon, -m-\mu-q- \varepsilon}(M,{\cal F},E)$,
we get
\begin{eqnarray*}
\|R_E(k)P& : &L^2(M,E)\rightarrow
H^{q+\varepsilon,p+\varepsilon_1}(M,{\cal F},E)\|\\
&\leq &\|R_E(k): H^{q+\varepsilon, -m-\mu-q- \varepsilon}(M,{\cal F},E) 
\rightarrow H^{q+\varepsilon,p+\varepsilon_1}
(M,{\cal F},E)\|\\
& & \|P: L^2(M,E)\rightarrow H^{-m,-\mu}(M,{\cal F},E)\|,
\end{eqnarray*}
that completes the proof.
\end{proof}

Let $A\in \Psi^m(M,E)$ be a transversally elliptic operator. We
assume that $A$ considered as an unbounded operator in $L^2(M,E)$ with
the domain $C^{\infty}(M,E)$ is essentially self-adjoint.
For any mesurable, bounded function $f$ on ${\Bbb R}$,
the bounded operator $f(A)$ in $L^2(M,E)$ is defined via the spectral theorem.
\begin{proposition}
\label{ftrace}
For any  mesurable, bounded function $f$ on ${\Bbb R}$
such that
$$
|f(\lambda)|\leq C(1+|\lambda|)^{-l}, \lambda\in {\Bbb R},
$$
with some $C>0$ and $l>q/m$, the operator
$f(A)$ is an operator of $G$-trace class with the following
estimate for its $G$-trace functional:
$$
|({\trG}(f(A)),k)|\leq C_{\varepsilon}
\|k\|_{ml,-ml,n+\varepsilon-ml}, k\in C^{\infty}_c(G).
$$
for any $\varepsilon>0$.
\end{proposition}

\begin{proof} 
By (\ref{decom}) and Remark~\ref{adding}, we have
\begin{equation}
\label{decom1}
A^{l}= P(l)+T(l),
\end{equation}
where $P(l)$ is an elliptic operator
of class $\Psi^{ml}(M,E)$ with the positive principal symbol,
and, for any real $s$, $T(l)$ defines a continuous mapping
$$
T(l):L^2(M,E)\rightarrow H^{s,-s-ml}(M,{\cal F},E).
$$
Let $Q\in\Psi^{-ml}(M,E)$ be a parametrix for $P$, that is, 
\begin{equation}
\label{par}
QP(l)=I-K, K\in \Psi^{-\infty}(M,E).
\end{equation}
Let $B\in \Psi^{0,-\infty}(M,{\cal F},E)$ and $u\in C^{\infty}(M,E)$. 
By (\ref{par}), we have the estimate
$$
\|Bf(A)u\|_{ml,t} \leq \|BQP(l)f(A)u\|_{ml,t} + \|BKf(A)u\|_{ml,t}
$$
The second term in the right-hand side of the last estimate can be estimated
as follows:
$$
\|BKf(A)u\|_{ml,t}\leq C\|B\|_{ml,t,t}\sup |f(\lambda)|\,\|u\|.
$$
Let us turn to the first term. By (\ref{decom1}), we have
$$
\|BQP(l)f(A)u\|_{ml,t}\leq \|BQA^{l}f(A)u\|_{ml,t}+
\|BQT(l)f(A)u\|_{ml,t}.
$$
Finally, the terms in the right-hand side of the last estimate can  
be estimated as follows:
\begin{eqnarray*}
\|BQA^{l}f(A)u\|_{ml,t}&\leq &C\|B\|_{ml,0,t}\sup |(1+|\lambda |)^{l}
f(\lambda)|\,\|u\|,\\
\|BQT(l)f(A)u\|_{ml,t}&\leq &C\|B\|_{ml,-ml,t}\sup |f(\lambda)|\,\|u\|.
\end{eqnarray*} 
Taking $t=n+\varepsilon-ml$ and applying Proposition~\ref{trace}, 
we complete the proof.
\end{proof}

\begin{remark}
Using Proposition~\ref{ftrace}, we can define distributional 
spectral invariants of transversally elliptic operators like as
a spectrum distribution function, a zeta function etc. We refer
the reader to \cite{trans} for analogous results for transversally
elliptic operators on manifolds equipped with a smooth action 
of a (noncompact) Lie group. 
\end{remark}

\subsection{Zeta-function}
\label{zetaS}
As in Section \ref{complex}, we assume that $A\in\Psi^{m}(M,E)$
is a transversally elliptic classical
pseudodifferential operator with the positive
transversal principal symbol, which is essentially self-adjoint,
invertible and positive definite in the Hilbert space $L^2(M,E)$
(see (T1) and (T2) above). By Proposition~\ref{ftrace}, for any
$\hbox{Re}\ z>q/m$, the operator $A^{-z}$ is an
operator of $G$-trace, and {\bf the distributional zeta-function}
of the operator $A$ is defined as follows:
\begin{equation}
\label{zeta}
\zeta_{A}(z)={\trG}(A^{-z}),\ \hbox{Re}\ z> q/m.
\end{equation}
Moreover, Proposition~\ref{ftrace} provides the following
estimate for the distributional zeta-function with any $\varepsilon>0$:
\begin{equation}
\label{zest0}
|(\zeta_{A}(z),k)|\leq C_{\varepsilon}
\|k\|_{q,-q,p+\varepsilon}, k\in C^{\infty}_c(G),
\hbox{Re}\ z> q/m.
\end{equation}

Now we turn to the problem of meromorphic continuation of
$\zeta_A(z)$. Actually, we consider a little bit more general situation.
\begin{theorem}
\label{zeta1}
Let $A$ as above and
$Q\in \Psi^{l,-\infty}(M,{\cal F},E)$, $l\in {\Bbb Z}$. Then the function
$z\mapsto \Tr(QA^{-z})$ is holomorphic for
${\rm Re}\, z> l+q/m$ and admits a (unique) meromorphic continuation
to ${\Bbb C}$ with at most simple
poles at points $z_k=k/m$ with integer $k\leq l+q$. Its residue at
$z=z_k$ is given by
$$
\res_{z=z_k}\ \tr (QA^{-z})=q \tau(QA^{-k/m}).
$$
\end{theorem}

\begin{proof} 
Using (\ref{decom}),
we construct a meromorphic continuation of
$\tr(QA^{-z})$ as follows:
$$
\tr(QA^{-z})= \TR(QP(z))+\tr(QT(z)).
$$
Here $QP(z)$ is a holomorphic family of operators of class
$\Psi^{mz+l,-\infty}(M,{\cal F},E)$ and the meromorphic extension
of its trace is given by $\TR(QP(z))$ due to Proposition~\ref{extension}.
Further, $QT(z)$ defines a continuous mapping
$QT(z):L^2(M,E)\rightarrow H^{s}(M,E)$
for any $s$. Therefore, by Proposition~\ref{trace}, the operator $QT(z)$
is of trace class for any $z\in {\Bbb C}$ with $\tr(QT(z))$, being an
entire function of $z$.
\end{proof}
As a corollary, we have the following result on a meromorphic continuation
of the zeta-function, considered as a distribution on $G$.
\begin{theorem}
\label{zeta0}
For any $k\in C^{\infty}_c(G)$, the zeta function
$(\zeta_A(z),k)$ of the operator $A$ extends to a
meromorphic function on the complex plane with simple poles at points
$z=q/m,(q-1)/m,\ldots.$
\end{theorem}
For future references,
we state the following theorem, which can be proved in the same manner
(see \cite[Proposition 4.2]{KV2}).
\begin{theorem}
\label{zeta2}
Let $A$ as above and
$Q_j\in \Psi^{l_j,-\infty}(M,{\cal F},E)$, $l_j\in {\Bbb Z}$, $j=1,\ldots,
N$. Then the function
$$
\zeta(z_1,\ldots,z_N)=\tr(Q_1A^{-z_1}\ldots Q_NA^{-z_N}), (z_1,\ldots,z_N)
\in {\Bbb C}^N
$$
admits a (unique) meromorphic continuation
to ${\Bbb C}^N$ with at most simple
poles on the hyperplanes $\sum_{j=1}^N l_j-m\sum_{j=1}^N z_j=k\in {\Bbb Z},
k\geq -q$.
Its residue at this hyperplane is given by
$$
\res \zeta(z_1,\ldots,z_N) =q \tau(Q_1A^{-z_1}\ldots Q_NA^{-z_N}).
$$
\end{theorem}

\section{Spectral triples of Riemannian foliations}
\subsection{Proof of main theorems}
\label{proof}
In this section, we complete proofs of our main theorems.
Recall \cite{Co-M,spview} that a spectral triple is a triple
$({\cal A}, {\cal H}, D)$, where:
\begin{enumerate}
\item ${\cal A}$ is an involutive algebra;
\item ${\cal H}$ is a Hilbert space equipped with
a $\ast$-representation of the algebra ${\cal A}$;
\item $D$ is a (unbounded) selfadjoint
operator in ${\cal H}$ such that
\medskip
\par
\noindent 1. for any $a\in {\cal A}$, the operator
$a(D-i)^{-1}$ is a compact operator in ${\cal H}$;
\medskip\par
\noindent 2. $D$ almost commutes with any $a\in {\cal A}$
in a sense that
$[D,a]$ is bounded for any $a\in {\cal A}$.
\end{enumerate}

One of the basic geometrical examples of spectral triples is given
by a triple $(A,{\cal H},D)$, associated with a compact Riemannian
manifold $M$:
\begin{enumerate}
\item The involutive algebra ${\cal A}$ is the algebra $C^{\infty}(M)$
of smooth functions on $M$;
\item The Hilbert space ${\cal H}$ is the $L^2$ space $L^2(M,\Lambda^*M)$
of differential forms on $M$, on which the algebra ${\cal A}$ acts
by multiplication;
\item The operator $D$ is the signature operator $d+d^*$.
\end{enumerate}

In this section, we will consider spectral triples
$({\cal A},{\cal H},D)$ associated with a compact foliated manifold
$(M,{\cal F})$:
\begin{enumerate}
\item The involutive
algebra ${\cal A}$ is the algebra $C^{\infty}_c(G)$;
\item The Hilbert space ${\cal H}$
is the space $L^2(M,E)$ of $L^2$-sections of a holonomy
equivariant Hermitian vector bundle $E$, on which
an element $k$ of the algebra ${\cal A}$ is represented via
the $\ast$-representation $R_E$;
\item The operator $D$ is a first order self-adjoint
transversally elliptic operator with the holonomy invariant transversal
principal symbol such that the operator $D^2$ is self-adjoint and
has the scalar principal symbol.
\end{enumerate}

The following theorem is Theorem~\ref{triple} of Introduction.

\begin{theorem}
\label{triple1}
Let $(M,{\cal F})$ be a closed foliated manifold.
Then a spectral triple $(A,{\cal H},D)$  as above
is a finite-dimensional spectral triple, that is:
\medskip\par
\noindent 1. for any $k\in C^{-\infty}_c(G)$, the operator
$R_E(k)(D-i)^{-1}$ is a compact operator in $L^2(M,E)$;
\medskip\par
\noindent 2. for any $k\in C^{-\infty}_c(G)$, $[D,R_E(k)]$ is
bounded in $L^2(M,E)$.
\end{theorem}

\begin{proof} 
The following proposition applied to $A=D-i$ implies,
in particular, the first part of this theorem.
\begin{proposition}
\label{invers}
Let $A\in \Psi^1(M,E)$ be a transversally elliptic operator,
invertible in the Hilbert space $L^2(M,E)$. Then,
for any $k\in {\rm OP}^{-n/q}({\cal F})$, the operator
$R_E(k)A^{-1}$ defines a continuous map from $L^2(M,E)$ to
$H^{1,p/q}(M,{\cal F},E)$. In particular, the operator
$R_E(k)A^{-1}$ is a compact operator in $L^2(M,E)$.
\end{proposition}
\begin{proof} By (\ref{param}), we have the following
representation:
\begin{equation}
\label{repr}
R_E(k)A^{-1}=R_E(k)P+R_E(k)R_2A^{-1},
\end{equation}
where $P\in  \Psi ^{-1}(M,E)$ and $R_2\in \Psi ^{0}(M,E)\bigcap
\Psi ^{-\infty}(N^*{\cal F},E)$. Since $k\in {\rm OP}^{-n/q}({\cal F})$,
by Proposition~\ref{13}, the first term, $R_E(k)P$, in the right-hand
side of (\ref{repr}) defines a continuous mapping from
$H^{s,t}(M,{\cal F},E)$ to $H^{s+1,t+n/q}(M,{\cal F},E)$ for any $s$ and $t$.
By (\ref{cont}), the operator $R_E(k)R_2A^{-1}$ defines a continuous
mapping from $L^2(M,E)$ to $H^{N,n/q-N}(M,{\cal F},E)$ for any $N$.
So we get that the operator $R_E(k)A^{-1}$
defines a continuous map from $L^2(M,E)$ to
$H^{1,p/q}(M,{\cal F},E)$.
\end{proof}

The second part of this theorem, concerning to boundedness of
commutators $[D,R_E(k)]$, follows from Proposition~\ref{com}.
\end{proof}

\begin{remark}
By Proposition~\ref{invers} and Proposition~\ref{trace},
it is easy to see that, for any $k\in C^{\infty}_c(G)$, the operator
$R_E(k)(D-i)^{-1}$ belongs to the Schatten ideal
${\cal L}^{q+\varepsilon}(L^2(M,E))$ for any $\varepsilon>0$,
therefore, the spectral triple
in question has the finite dimension $d$, which is equal to $q$.
\end{remark}

Now we turn to a description of the dimension spectrum for the
spectral triples under consideration.
First, recall briefly the definition of the
dimension spectrum \cite{Co-M,spview}.
Let $({\cal A}, {\cal H}, D)$ be a spectral triple.
Denote by $\delta$ an (unbounded) derivation on the algebra
${\cal L}({\cal H})$
of all linear operators in ${\cal H}$, given by the formula
\begin{equation}
\label{derivative}
\delta(T)=[|D|,T], T\in {\cal L}({\cal H}).
\end{equation}
Assume that, for any $a\in {\cal A}$,
\begin{equation}
\label{restr}
a\in \bigcap_{n>0} {\rm Dom}\, \delta^n,
[D,a]\in \bigcap_{n>0} {\rm Dom}\, \delta^n,
\end{equation}
and denote by ${\cal B}$ the algebra generated
by the elements $\delta^n(a), a\in{\cal A}, n\in {\Bbb N}$.
Then the operator $b|D|^{-z}$ is of trace
class for $Re\ z>d$, where $d$ is the top spectrum dimension and $b\in
{\cal B}$, and we can define the distributional zeta function $\zeta_b(z)$ of
the operator $|D|$ by the formula
$$
\zeta_b(z)=\tr(b|D|^{-z}), b\in {\cal B}, {\rm Re}\, z>d.
$$

\begin{definition}  
A spectral triple $({\cal A},{\cal H},D)$ has {\bf discrete
dimension spectrum} ${\rm Sd}\subset{\Bbb C}$, if ${\rm Sd}$
is a discrete subset in ${\Bbb C}$, the triple satisfies the 
assumptions~(\ref{restr}), and, for any $b\in {\cal B}$, the
distributional zeta function $\zeta_b(z)$ extends
holomorphically to ${\Bbb C}\backslash {\rm Sd}$ such that $\Gamma(z)
\zeta_b(z)$ is of rapid decay on vertical lines $z=s+it$,
for any $s$ with ${\rm Re}\,s >0$.

The dimension spectrum is said to be {\bf simple}, if the singularities
of the function $\zeta_b(z)$ at $z\in {\rm Sd}$ are at most simple poles.
\end{definition}

From now on, let $D$ be a first order transversally elliptic
operator with the holonomy invariant transversal
principal symbol such that the
operator $D^2$ has the scalar principal symbol, and $D$ and $D^2$ are
self-adjoint.
We start with a description of
the domain of the derivative
$\delta : C^{\infty}_c(G)\rightarrow {\cal L}(L^2(M,E))$, given
by (\ref{derivative}).
To do this, we slightly refine the description of the operator
$|D|$, given by Proposition~\ref{powers}.

\begin{lemma}
\label{powers1}
Under the current assumptions, the operator $|D|$ has the form
\begin{equation}
\label{mod}
|D|= P+T,
\end{equation}
where $P$ is a pseudodifferential operator
of class $\Psi^{1}(M,E)$ with the scalar principal symbol and
the holonomy invariant transversal principal symbol,
and $T$ defines a continuous mapping
\begin{equation}
\label{map}
T:H^{s,k}(M,{\cal F},E)\rightarrow H^{s+t,k-t-1}(M,{\cal F},E)
\end{equation}
for any $s, k$ and $t$.
\end{lemma}
\begin{proof} By Proposition~\ref{powers}, the operator
$|D|=(D^2)^{1/2}$ can be represented in the form (\ref{mod}) with
$P$ and $T$, satisfying almost all the conditions stated in the lemma.
It only remains to prove (\ref{map}) for any $s, k$ and $t$, because,
by Proposition~\ref{powers}, we know it is true
for any $s,k$ with $s+k=0$. Recall that $T$ is given by
$$
T = \frac{i}{2\pi} \int_{\Gamma} \lambda ^{-1/2}D^2T(\lambda )d\lambda ,
$$
and $T(\lambda)$ is given by (\ref{Tlambda}), therefore, it is easy to see
that, in order to prove (\ref{map}), it suffices to state the standard
resolvent estimate in any Sobolev space:
\begin{equation}
\label{rez:Sob}
\|(D^2-\lambda)^{-1}:H^s(M,E)\rightarrow H^s(M,E)\|\leq C/|\lambda|
\end{equation}
for any real $s$ and $\lambda\in\Lambda_{\delta}$ with $|\lambda|$ large
enough ($\delta\in (0,\pi)$ is arbitrary).

Recall that the estimate (\ref{rez:Sob}) for $s=0$ is a direct consequence
of self-adjointness of $D^2$. Let $\Lambda_s=(I+\Delta_M)^{s/2}$.
Then $\Lambda_sD\Lambda_{-s}=
D+B_s$, where $B_s$ is a bounded operator in $L^2(M,E)$, and
$\Lambda_sD^2\Lambda_{-s}=D^2+DB_s+B_sD+B_s^2$. It can be easily checked
that, for any $\varepsilon>0$, we have the estimate
$$
((DB_s+B_sD+B_s^2)u,u)\leq \varepsilon\|Du\|^2+C_{\varepsilon}\|u\|^2,
u\in C^\infty(M,E),
$$
which implies the estimate (\ref{rez:Sob}) for any $s$ due to well-known
facts of the perturbation theory of linear operators (see, for instance,
\cite{Kato}).
\end{proof}

\begin{proposition}
\label{dom}
Under the current assumptions,
the operator $\delta(K)=[|D|,K]$ is an operator of
class $\Psi^{0,-\infty}(M,{\cal F},E)$
for any $K\in \Psi^{0,-\infty}(M,{\cal F},E)$.
\end{proposition}

\begin{proof}
By (\ref{mod}), $\delta(K)=[P,K]+[T,K]$,
where $P$ and $T$ are as in Lemma~\ref{powers1}.
Since $P$ has the scalar principal symbol and the holonomy invariant
transversal principal symbol, Proposition~\ref{module} implies that the
operator $[P,K]$ belongs to the class
$\Psi^{0,-\infty}(M,{\cal F},E)$. By (\ref{map}), the operator
$[T,K]$ is a smoothing operator in the scale $H^{s,k}(M,{\cal F},E)$,
and, therefore, belongs to $\Psi^{-\infty}(M,E)$.
\end{proof}
By Proposition~\ref{dom}, we get immediately the following
characterization of the domain of the derivative $\delta$.

\begin{proposition}
\label{dom1}
Given an operator $D$ as above, the class
$\Psi^{0,-\infty}(M,{\cal F},E)+{\rm OP}^{-1}({\cal F})$ is contained
in the domain of the derivative $\delta$.
\end{proposition}

For further references, we note the following, slightly more general
assertion of the domain of the derivative $\delta$, which follows
from the estimates of Proposition~\ref{com}.

\begin{proposition}
\label{dom2}
Given a first order transversally elliptic operator $D$ with the
holonomy invariant transversal principal symbol,
the class ${\rm OP}^{-1}({\cal F})$ is contained in the domain of $\delta$.
\end{proposition}

Now we are ready to prove the second main result of the paper,
Theorem~\ref{dim} of Introduction.

\begin{theorem}
\label{dim1}
A spectral triple $({\cal A},{\cal H},D)$ as in Theorem~\ref{triple1}
has discrete spectrum dimension ${\rm Sd}$, which is contained in the set
$\{v\in {\Bbb N}:v\leq q\}$ and is simple.
\end{theorem}

\begin{proof}
First of all, we have to verify (\ref{restr}).
Given ${\cal A}$, being the algebra $C^{\infty}_c(G)$,
by Propositions~\ref{dom} and \ref{dom1},
all the assumptions on ${\cal A}$ formulated in (\ref{restr}) are satisfied,
that is, for any $k\in C^{\infty}_c(G)$,
$R_E(k) \in \bigcap_{n>0} {\rm Dom}\, \delta^n,
[D,R_E(k)] \in \Psi^{0,-\infty}(M,{\cal F},E)
\subset \bigcap_{n>0} {\rm Dom}\, \delta^n$.
Moreover, the algebra ${\cal B}$,
generated by the elements $\delta^n(R_E(k)), k\in{\cal A}, n\in {\Bbb N}$,
is contained in $\Psi^{0,-\infty}(M,{\cal F},E)$,
The rest of the proof follows immediately from Theorem~\ref{zeta1}.
\end{proof}

\subsection{Geometric example}
\label{self}
In this section, we discuss an example of a spectral triple given by
the transverse signature operator on a Riemannian foliation.

Let $(M,{\cal F})$ be a Riemannian foliation equipped with a
bundle-like metric $g_M$. Let $F=T{\cal F}$ be the tangent bundle
to ${\cal F}$, and $H=F^{\bot}$ be the orthogonal
complement to $F$. So we have a decomposition of
$TM$ into a direct sum
\begin{equation}
\label{decomp1}
TM=F\oplus H.
\end{equation}
The de Rham differential $d$ inherits the
decomposition~(\ref{decomp1}) in the form
$$
d=d_F+d_H+\theta.
$$
Here the tangential de Rham differential $d_F$
and the transversal de Rham differential $d_H$ are first order differential
operators, and $\theta$ is zeroth order. Moreover, the
operator $d_F$ doesn't depend on a choice of $g_M$
(see, for instance, \cite{Re}).

The conormal bundle $N^*{\cal F}$ has a leafwise flat connection (the Bott
connection) defined by the lifted foliation ${\cal F}_N$. The parallel
transport along leafwise
paths with respect to this connection defines a representation of the
holonomy groupoid $G$ in fibres of $N^*{\cal F}$.
Since $H^*\cong N^*{\cal F}$, the bundles $H^{*}$ and $\Lambda^{*}H^{*}$
are holonomy equivariant.

Define a triple $({\cal A}, {\cal H}, D)$ to be given by
the space ${\cal H}=L^2(M,\Lambda^{*}H^{*})$ of the transversal
differential forms equipped with the action of the
algebra ${\cal A}=C^{\infty}_c(G)$ and by the transverse signature
operator $D=d_H+d^*_H$.

Let us check that this spectral triple satisfy all the assumptions of Section
\ref{proof}. Under the isomorphism $H^*\cong N^*{\cal F}$,
the transversal principal symbol, $\sigma_D(\eta)\in
{\cal L}(\Lambda^*N^*{\cal F}$, of $D$ is given by the formula
$$
\sigma_D(\eta)=e_{\eta}+i_{\eta},\ \eta\in \tilde{N}^*{\cal F}
$$
($e_\eta$ and $i_\eta$ are the exterior and the interior multiplications by
$\eta$ accordingly),
from where one can easily see holonomy invariance of $\sigma_D$.
We also have
$$
\sigma_{D^2}(\eta)=|\eta|^2I_{\pi(\eta)},\ \eta\in \tilde{N}^*{\cal F}.
$$
Finally, essential self-adjointness of $D$ and $D^2$ follows from
the finite propagation speed arguments of \cite{Chernov}.

\subsection{Concluding remarks}
1. Using the well-known relationship between the zeta-function and the
heat trace, we can derive from Theorem~\ref{zeta0} the following fact on
heat trace asymptotics for transversally elliptic operators, which was
used in \cite{asymp}.

\begin{proposition}
\label{Proposition 4.3}
Let $(M,\cal F)$ be a compact foliated manifold and $E$ be
an Hermitian vector bundle on $M$.
Given a differential operator $P\in  \Psi ^{m}(M,E), m>0$,
with the positive transversal principal symbol, which is self-adjoint and
positive in $L^2(M,E)$, we have an asymptotic expansion of the $G$-trace of
the parabolic semigroup operator $e^{- tP}$:
$$
({\trG}(e^{- tP}),k)=\tr R_E(k)e^{- tP},
k\in C^{\infty}_{c}(G),
$$
given by the formula
$$
{\trG}(e^{- tP})\sim
\sum_{l=0}^{\infty} a_{l}t^{-q/m+l/m}, t\rightarrow 0+,
$$
where $a_l,l=0,1,\ldots$ are distributions on $G$ with $a_0$,
given by the formula:
\begin{equation}
\label{lead}
(a_0,k) = \int_{M} (\int_{N_{x}^{\ast}{\cal F}}\
\Tr e^{-\sigma _{P}(\eta )} k(x)d\nu (\eta ))dx,
k\in C^{\infty}_{c}(G).
\end{equation}
\end{proposition}

2. There are two observations, based on the explicit formula
(\ref{lead}) for the leading coefficient in the heat trace expansion, or,
equivalently, for the residue of the distributional zeta-function at
the point $z=-q/m$, which have an interesting interpretation in terms of
noncommutative spectral geometry.

Assume that $(M,{\cal F})$ is a compact Riemannian foliated manifold
with a bundle-like metric $g_M$, and a spectral triple $({\cal A},
{\cal H}, D)$ be given by the transverse signature operator (as in Section
\ref{self}).

Denote by  $C^{*}_E(G)$ the closure of the
set $R_E(C^{\infty}_c(G))$ in the uniform operator topology
of ${\cal L}(L^2(M,E))$.
It is easy to see that the $\ast$-homomorphism
$R_E: C^{\infty}_c(G)\rightarrow C^{*}_E(G)$
extends to a $\ast$-homomorphism of the full
$C^{*}$-algebra $C^{*}(G)$, 
$R_E: C^{*}(G)\rightarrow C^{*}_E(G)$.
By \cite{F-Sk}, we have the natural projection
$\pi_E : C^{\ast}_{E}(G) \rightarrow  C^{\ast}_{r}(G)$.
If ${\cal I}$ is an involutive ideal in $C^{*}_E(G)$, we may think
of the spectral triple $({\cal I}, {\cal H}, D)$ as a subset of our
spectrally defined geometrical space and look for its dimension spectrum.
Then Theorem~\ref{dim} implies the following fact about the dimension spectra
of "subsets" in the singular space $M/{\cal F}$.

\begin{proposition}
Let ${\cal I}$ be an involutive ideal in $C^{*}_E(G)$.
Then $q\in {\rm Sd}({\cal I})$ iff $\pi_E({\cal I})\not=0$.
In particular, if
$\pi_E({\cal I})=0$, then the top spectral dimension
of the spectral triple $({\cal I}, {\cal H}, D)$ is less than  $q$.
\end{proposition}
This fact can be also interpreted as a fact about a noncommutative analogue
of the integral in the case under consideration.
Let $I$ be a functional
on $C^{\infty}_c(G)$, given by the formula
$$
I(k)=\tau(R_E(k)|D|^{-q}), k\in C^{\infty}_c(G).
$$

\begin{proposition}
Under the current assumptions, the functional
$I$ is given by the following formula
$$
I(k)
=\frac{q}{\Gamma(\frac{q}{2}+1)}{\trF}\pi_E(k), k\in C^{\infty}_c(G),
$$
and can be extended by continuity to a functional on $C^{*}$-algebra
 $C^{*}_E(G)$.
In particular, $I(k)=0$ for any  $k\in C^{*}_E(G), \pi_E(k)=0$.
\end{proposition}

Otherwise speaking, the functional
$I$ coincides on $C^{\infty}_c(G)$ (up to some multiple) with the von Neumann
trace ${\trF}$, given by the Riemannian transversal volume due to
the noncommutative integration theory
\cite{Co79}, and the support of $I$ is a "regular" part of our geometrical
space.

\begin{acknowledgement}
The work was done during visits to the Max Planck Institut
f\"{u}r Mathematik at Bonn and to Institut des Hautes Etudes Scientifiques.
I wish to express my gratitude to these institutes for their hospitality 
and support.
\end{acknowledgement}


\begin{thebibliography}{99}
\bibitem{A-U} Antoniano, J., Uhlmann, G.: A functional calculus for a
class of pseudodifferential operators with singular symbols.
Proc. Symp. Pure Math. {\bf 43}, 5-16(1985).
\bibitem{Atiyah} Atiyah, M.F.: Elliptic operators and compact groups.
(Lecture Notes in Math. 401) Berlin Heidelberg New York: Springer 1974.
\bibitem{B-G} Beals, R., Greiner, P.: Calculus on Heisenberg 
manifolds.
(Annals of Math. Studies 119) Princeton: Princeton Univ. Press 1988.
\bibitem{Chernov} Chernoff, P.R.: Essentiall self-adjointness of powers
of generators of hyperbolic equations. J. Funct. Anal. {\bf 12}(1973),
401-414.
\bibitem{Co79} Connes, A.: Sur la theorie  non  commutative  de  
l'integration. In: Lecture Notes in Math. 725. Berlin Heidelberg New York: 
Springer 1979. - 19 - 143.
\bibitem{CoNG} Connes, A.: Noncommutative differential geometry. Inst.
Hautes Etudes Sci. Publ. Math. {\bf 62}, 257-360(1985).
\bibitem{Co} Connes, A. Noncommutative geometry. London: Academic Press 1994.
\bibitem{spview} Connes, A.: Geometry from the spectral point of
view. Lett. Math. Phys. {\bf 34}, 203-238(1995). 
\bibitem{Co-M} Connes, A., Moscovici, H.: The local index
formula in noncommutative geometry. Geom. and Funct. Anal. {\bf 5}, 
174-243(1995).
\bibitem{F-Sk} Fack, T., Skandalis, G.: Sur les representations 
et ideaux
de la $C^*$-algebre d'un feuilletage. J. Operator Theory {\bf 8}, 
95-129(1982).
\bibitem{GS79} Guillemin, V., Sternberg, S.: Some problems in integral
geometry and some related problems in microlocal analysis. Amer. J. Math.
{\bf 101}, 915-959(1979).
\bibitem{Gu_1} Guillemin, V.: Gauged Lagrangian distributions. 
Adv. Math. {\bf 102}, 184-201(1993).
\bibitem{Gu85} Guillemin, V.: A new proof of Weyl's formula on the
asymptotic distribution of eigenvalues. Adv. Math. {\bf 55}, 131-160(1985).
\bibitem{Gu_2} Guillemin, V.: Residue traces for certain algebras 
of Fourier
integral operators. J. Funct. Anal. {\bf 115}, 391-417(1993).
\bibitem{Gu-U} Guillemin, V., Uhlmann, G.: Oscillatory integrals 
with singular
symbols. Duke Math. J. {\bf 48}, 251-267(1981).
\bibitem{Gr-U} Greenleaf, A., Uhlmann, G.: Estimates for singular Radon
transforms and pseudodifferential operators with singular symbols. J. Funct.
Anal. {\bf 89}, 202-232(1990).
\bibitem{H-Sk} Hilsum, M., Skandalis, G.: Morphismes $K$-orientes 
d'espaces de
feuilles et fonctorialite en theorie de Kasparov. Ann. Sci. Ecole Norm. Sup.
(4) {\bf 20}, 325-390(1987).
\bibitem{Ho} H\"ormander, L.: The analysis of linear
partial differential operators I. Berlin Heidelberg
New York Tokyo: Springer 1983.
\bibitem{H3} H\"ormander, L.: The analysis of linear
partial differential operators III. Berlin Heidelberg
New York Tokyo: Springer 1985.
\bibitem{H4} H\"ormander, L.: The analysis of linear
partial differential operators IV. Berlin Heidelberg
New York Tokyo: Springer 1986.
\bibitem{Kato} Kato, T.: Perturbation theory for linear
operators. Berlin: Springer, 1980.
\bibitem{KV1} Kontsevich, M., Vishik, S.: Determinants
of elliptic pseudo-differential operators.
Preprint MPI/94-30, 1994, 156pp.
\bibitem{KV2} Kontsevich, M., Vishik, S.: Geometry of  determinants 
of elliptic operators. In: Functional analysis on the eve of the 21st
century. Vol. I. (Progress in Mathematics. Vol. 132) Boston:
Birkh\"auser 1996. - 173 - 197.
\bibitem{trans}Kordyukov, Yu.A.: Transversally elliptic
operators on $G$-manifolds of bounded geometry.
Russ. J. Math. Ph. {\bf 2}, 175-198(1994); {\bf 3}, 41-64(1995).
\bibitem{tang} Kordyukov, Yu.A.:
Functional calculus for tangentially elliptic operators on
foliated manifolds, In: Analysis and Geometry in Foliated
Manifolds, Proceedings of the VII International Colloquium on Differential
Geometry, Santiago de Compostela, 1994. - Singapore: World Scientific 1995.-
113-136.
\bibitem{asymp} Kordyukov, Yu.A.: Semiclassical spectral
asymptotics on a foliated manifold, to appear.
\bibitem{noncom2} Kordyukov, Yu.A.: Noncommutative spectral 
geometry of
transversally triangular Riemannian foliations, in preparation.
\bibitem{Molino} Molino, P.: Riemannian Foliations, Progress in Math.
no. 73, Basel Boston: Birkh\"auser, 1988.
\bibitem{Re} Reinhart, B.L.: Differential Geometry of
Foliations. Berlin Heidelberg New York: Springer 1983.
\bibitem{Wo} Wodzicki, M.: Noncommutative residue. Part I. Fundamentals.
In: K-theory, arithmetic and geometry (Moscow, 1984-86), 
Lecture Notes in Math. 1289. Berlin Heidelberg New York: Springer 1987.
-320-399.
\end{thebibliography}
\end{document}